\documentclass[1p, final]{elsarticle} %[review]

\usepackage{lineno,hyperref}
\modulolinenumbers[5]

\usepackage[cmex10]{amsmath}
\usepackage{amssymb}
\usepackage{bm}
\usepackage{algorithmicx}
\usepackage{algpseudocode}
\usepackage{algorithm}

\usepackage[caption=false,font=normalsize,labelfont=sf,textfont=sf]{subfig}

\usepackage{diagbox}
\usepackage{multirow}
\usepackage{color}

\hypersetup{
    colorlinks=true,       % false: boxed links; true: colored links
    linkcolor=black,          % color of internal links (change box color with linkbordercolor)
    citecolor=black,        % color of links to bibliography
    filecolor=black,      % color of file links
    urlcolor=black           % color of external links
}

\begin{document}

\begin{frontmatter}

\title{FLORIS and CLORIS: Hybrid Source and Network Localization Based on Ranges and Video}

%\author{{Beatriz~Quintino~Ferreira}\corref{mycorrespondingauthor}, Jo\~{a}o~Gomes, Cl\'{a}udia Soares, Jo\~{a}o~P.~Costeira}

\author{Beatriz~Quintino~Ferreira\corref{mycorrespondingauthor}}
\cortext[mycorrespondingauthor]{Corresponding author}
\ead{beatrizquintino@isr.tecnico.ulisboa.pt}
\author{Jo\~{a}o~Gomes}
\author{Cl\'{a}udia Soares}
\author{Jo\~{a}o~P.~Costeira}

\address{Institute for Systems and Robotics -- LARSyS, Instituto Superior T\'{e}cnico, Universidade de Lisboa, Portugal.}

\begin{abstract}
We propose hybrid methods for localization in wireless sensor networks fusing noisy range measurements with angular information (extracted from video). Compared with conventional methods that rely on a single sensed variable, this may pave the way for improved localization accuracy and robustness. We address both the single-source and network (i.e., cooperative multiple-source) localization paradigms, solving them via optimization of a convex surrogate. The formulations for hybrid localization are unified in the sense that we propose a single nonlinear least-squares cost function, fusing both angular and range measurements. We then relax the problem to obtain an estimate of the optimal positions. This contrasts with other hybrid approaches that alternate the execution of localization algorithms for each type of measurement separately, to progressively refine the position estimates.
Single-source localization uses a semidefinite relaxation to obtain a one-shot matrix solution from which the source position is derived through factorization. Network localization uses a different approach where sensor coordinates are retained as optimization variables, and the relaxed cost function is efficiently minimized using fast iterations based on Nesterov's optimal method.
Further, an automated calibration procedure is developed to express range and angular information, obtained through different devices, possibly deployed at different locations, in a single consistent coordinate system. This drastically reduces the need for manual calibration that would otherwise negatively impact the practical usability of hybrid range/video localization systems.
We develop and test, both in simulation and experimentally, the new hybrid localization algorithms, which not only overcome the limitations of previous fusing approaches, but also compare favourably to state-of-the-art methods, even outperforming them in some scenarios.
\end{abstract}


\end{frontmatter}

%\linenumbers

\section{Introduction}\label{sec:introduction}
The ``where am I'' problem has always been a key issue in the field of technology, both for human mobility and for robots/autonomous vehicles. Nonetheless, there are several scenarios, such as indoors or underwater, in which the most popular localization system, the Global Positioning System (GPS), is not available and where location awareness will soon become an essential feature. These environments pose challenges such as strong multi-path/non line-of-sight propagation, diffractions or interferences, which lead to over-meter accuracy for the majority of existing systems. Such accuracy may be insufficient for numerous applications, and the key to overcoming this issue may lie on exploring hybrid schemes~\cite{You:are:here:SPECTRUM}.

Focusing on indoor environments, most of the proposed localization systems use only one type of measurement, typically range~\cite{patwari:2005}. Yet, as wireless sensor networks (WSN) are becoming ubiquitous, it makes sense to try to infer positions from the spatial cues provided by the various sensors on-board networked devices to improve the accuracy and/or coverage. In this vein, the methods introduced in our work address centralized localization problems fusing distances (obtained acoustically or via radio signals) and angular information (emphasizing the use of video cameras).

Our approach is based on convex optimization and relaxation techniques, providing a sound framework for dealing with potentially very noisy measurements from low-cost network nodes. We propose a unified framework whereby range and angular information is incorporated into a single nonlinear least-squares (LS) nonconvex cost function, and the position estimates are obtained by relaxing it and finding the global optimum using a single minimization procedure. This is preferable to hybrid ``ping-pong'' iterative refinement schemes that alternate between localizing sources/sensors using a single type of measurement. These ad-hoc schemes may oscillate over time and require that an initial configuration be estimated from one type of measurement alone, whereas in some configurations of interest the number of available ranges or bearings, taken independently, may be insufficient to determine the position unambiguously.

We highlight the following contributions of our work:
\begin{itemize}
\item A hybrid single-source localization method (FLORIS --- Fused LOcalization using Ranges and Incident Streaks) where the original optimization problem is relaxed to a semidefinite program (SDP) whose one-shot\footnote{The term one-shot is used in the sense that only one problem is sent to the solver.} matrix solution may be calculated by a general-purpose convex solver. The relaxation is tight, yielding a high-accuracy estimate for the source position through matrix factorization. A description and preliminary performance characterization of this algorithm were presented in~\cite{Ferreira:hybrid_local:2015};
\item A hybrid network localization method (CLORIS --- Cooperative LOcalization using Ranges and Incident Streaks) based on a so-called disk relaxation of the nonconvex and nonlinear LS cost function. This approach retains the sensor coordinates as optimization variables, but explicitly requires iterative processing to converge to the global optimum of the convexified cost function. Its gradient-based nature allows a parallel implementation, and Nesterov's optimal method leads to a very fast and accurate algorithm, whose iterations converge quickly with similar complexity to much less efficient gradient descent methods;
\item An automated self-calibration procedure that expresses range and angular information, obtained through ranging sensors and video cameras, possibly deployed at different locations, in a single consistent coordinate system. This is a critical enabling component for our unified approach, where the cost functions include both range and angular terms expressed in a common frame. Manual calibration of the sensor networks is somewhat cumbersome, and streamlining the procedure is essential to ensure the practical usability of a hybrid range/video localization system.
\end{itemize} 

The new methods, which operate in 2D, 3D or even higher arbitrary dimensions, were fully tested both in simulation and in real experiments and achieved very encouraging results. In particular, they outperformed other benchmark methods when measurements were quite noisy.

Throughout, both scalars and individual position vectors will be represented by lower-case letters. Vectors of concatenated coordinates and matrices will be denoted by boldface lower-case and upper-case, respectively. The superscript $(\cdot)^*$ stands for conjugate transpose and $(\cdot)^T$ for the transpose of the given vector or matrix. $\otimes$ represents the Kronecker product, and $\|\mathbf{A}\|_F$ the Frobenius norm of matrix $\mathbf{A}$. $\mathbf{I}_m$ is the identity matrix of size $m \times m$ and $\mathbf{1}_m$ is the (column) vector of $m$ ones. For symmetric matrix $\mathbf{X}$, $\mathbf{X} \succeq 0$ means that $\mathbf{X}$ is positive semidefinite. The cardinality of set $\mathcal{A}$ is denoted by $|\mathcal{A}|$.

\subsection{Related Work}
In WSN localization, range information can be measured from absolute or differential travel times~\cite{sensnet:bach:2005}, and usually produces robust results for ranges up to about 10 meters. It can also be inferred, much less reliably, from received signal power~\cite{rss:pinar:2013}. Orientation (Angle of Arrival) is used less frequently for localization in WSN~\cite{Biswas05integrationof}, but remains a relatively popular alternative for outdoor geolocation and navigation when GPS is unavailable or unreliable~\cite{Ho:Bearing_only:2008}. In indoor scenarios orientation is a key enabler for augmented-reality systems, which superimpose realistic synthetic objects on camera images. 
%The widely popular ARUCO library, e.g., provides black-box functionality to extract orientation from images with the aid of simple fiducial markers~\cite{Aruco2014}, and is used in our work as well.
Distance and orientation retrieved from video are very precise and reliable at short ranges and in the absence of occlusion~\cite{cameraphones:mulloni:2009, alahi:2015:rgb-w}, whereas ranging devices such as the acoustic Cricket system used in our setup have moderate precision over an extended operating range~\cite{cricketlocalsys:Priyantha:2000}. The complementary strengths of these sensors make them extremely appealing to be used in synergy and seamlessly, paving the way to more accurate localization solutions.

Iterative and ad-hoc methods for WSN localization that combine range and angle information have been proposed in~\cite{Urruela:AOATDOAAdHoc:2006,Huang:CRLB_AOA:2016, Jia:TDOA_AOA_2018, Park:RangeAngleAd-Hoc:2013,Kleunen:TOFDOAunderwater:2014} (the latter being specific for single-source localization in underwater scenarios).
In contrast with the previous iterative schemes, and to the best of our knowledge, only a few recent attempts were made to genuinely fuse range and orientation for hybrid localization~\cite{Tomic:RSS_AOA:2016,Tomic:DAT:2017,Biswas05integrationof,alessio:selfcalib:2012}. The first two of these bases localization on received signal strength (RSSI) and angle measurements, whereas the last two are closer to our approach. However, these methods impose severe limitations; the one in~\cite{Biswas05integrationof} is SDP-based, but the way it encodes angular constraints is specific for 2D and very different from the techniques used in FLORIS and CLORIS. The approach of~\cite{alessio:selfcalib:2012} uses a bilinear matrix formulation inspired in computer vision techniques, but the assumption that the range and visual anchors overlap is too restrictive for many WSN scenarios of interest. %scenarios
In this work, we succeed in overcoming the limitations of both~\cite{Biswas05integrationof} and~\cite{alessio:selfcalib:2012} by deriving a novel formulation based on a single optimization problem that jointly accounts for range and bearing data obtained from \emph{arbitrarily placed} heterogeneous sensors in \emph{2D or 3D}, the latter being particularly relevant in indoor scenarios. Optimal placement of nodes - either static~\cite{sensnet:bach:2005,Biswas05integrationof,Aspnes:TheoryNetLoc:2006} or mobile~\cite{Lin:DeploymentHetero:2016} - can significantly enhance the estimation/localization accuracy, but is beyond the scope of our work.  

Semidefinite programming and relaxation techniques~\cite{convexopt:boyd:2004} have been successfully used before in range-based localization to obtain high-quality approximations to the maximum of the nonconvex likelihood function under Gaussian noise (see~\cite{cheung:2004:icassp, SLCP:pinar:2014} and references therein). The approach taken here for the hybrid single-source case builds upon the Source Localization with Nuclear Norm (SLNN) algorithm for range-only measurements~\cite{SLCP:pinar:2014}, which provides state-of-the-art performance among SDP methods. The key idea is to reformulate the Maximum-Likelihood (ML) problem as finding a suitably constrained set of directions such that they nearly intersect (at the estimated source position). While we follow the general approach of SLNN, we adopt specific reformulations and relaxations for the modified LS cost function that includes additional terms to account for angle measurements. Importantly, SLNN's high accuracy and robustness to noise carries over to the proposed method.

We also address the sensor network paradigm, in which multiple network nodes (sensors) do not know their positions, as opposed to anchor nodes whose positions are known \emph{a priori}. This paradigm is very compelling and has emerged in numerous applications since it enhances accuracy and robustness of localization while circumventing the need for high-density anchor placement and high power allocation~\cite{patwari:2005,Conti:NetworkExperimentation:2012}. Similarly to the single-source case, centralized range-based methods that estimate all network positions based on the minimization of a LS (i.e., Gaussian ML) criterion and also on its weighted versions were proposed in~\cite{colabLS:shang:2004, colabWLS:destino:2011}. Semidefinite relaxations for the original nonconvex ML problem have also been developed for line-of-sight and non-line-of-sight environments~\cite{EDM_robust:pinar:2011,Vaghefi:CooperatoiveNLOS:2016}. However, the resulting SDPs may be large, limiting the problem sizes that can reasonably be tackled, at present, to fewer than a hundred nodes. With the advent of large sensor networks this becomes too restrictive, and alternative approaches, often designed for distributed processing between nodes in the network, have been proposed. In~\cite{Hero:POCS:2005} distributed single-source localization is formulated as a convex feasibility problem and solved via projections onto convex sets. This method was later extended to the problem of sensor network localization in~\cite{Gholami:coop_implicit_convex:2013}. Another parallel approach was proposed in~\cite{colab:dist:parallel:2014} based on two consecutive relaxations of the ML function (one SDP, followed by an edge-based relaxation). A simpler algorithm, based on the Gauss Seidel framework and possessing convergence guarantees, was proposed in~\cite{colab:sequen:2010}. However, the algorithm operates sequentially across nodes, thus requiring the existence of a global coordination mechanism for in-network processing. Refinement methods for nonconvex costs such as~\cite{Xia:dist_adaptive_pos:2016} require good initial position estimates to converge and are not reviewed here.

 % The distributed state-of-the-art methods can be further divided into two subsets: initialization dependent and initialization independent. The first set comprises methods which directly solve the original nonconvex problem, as~\cite{colab:initdepend:2010}. In the latter, there is an intermediate convex relaxation step, with a solution that must approximate the global problem solution, regardless the initialization given. Such initialization independent approach includes both sequential, as in~\cite{colab:sequen:2010}, and parallel methods, as in~\cite{colab:dist:parallel:2014} and~\cite{soares:coopdistlocal:2014}, showing convergence guarantees.

The parallel method proposed in~\cite{soares:coopdistlocal:2014} outperforms previous cooperative localization methods~\cite{EDM_robust:pinar:2011, colab:dist:parallel:2014, Gholami:coop_implicit_convex:2013} not just in accuracy but also in computational and communication efficiency. This fast and scalable method minimizes the convex underestimator of the ML cost by resorting to a simple relaxation of terms that replaces circles (nonconvex sets) with disks, while retaining the same variables (sensor positions) as in the original ML formulation. The new hybrid network localization approach presented here stems from the range-only algorithm in~\cite{soares:coopdistlocal:2014}, introducing a modification relative to the angle measurements available between pairs of sensors, or between sensors and anchors. The new method inherits the accuracy, efficiency and scalability of the original range-based method of~\cite{soares:coopdistlocal:2014}. CLORIS also remains a fully parallel algorithm that is suitable for in-network processing.

\subsection{Problem Formulation}

In the single-source case, let $x \in \mathbb{R}^n$ denote the unknown position of the source, to be estimated based on measurements to a set of $m$ reference points (anchors) $a_k \in \mathbb{R}^n$, $k = 1, \, \ldots, m$, whose positions and orientations are known. Of these, the ones whose indices belong to set $\cal{R}$ provide range measurements to the source, $d_k = \|x - a_k \|$, whereas those with indices in $\cal{T}$ measure bearings encoded by unit-norm vectors, $u_j = \frac{x - a_j}{\|x - a_j\|}$. Actual measurements are corrupted by noise.

\begin{figure}[h]
  \centering
  \includegraphics[width=0.42\columnwidth]{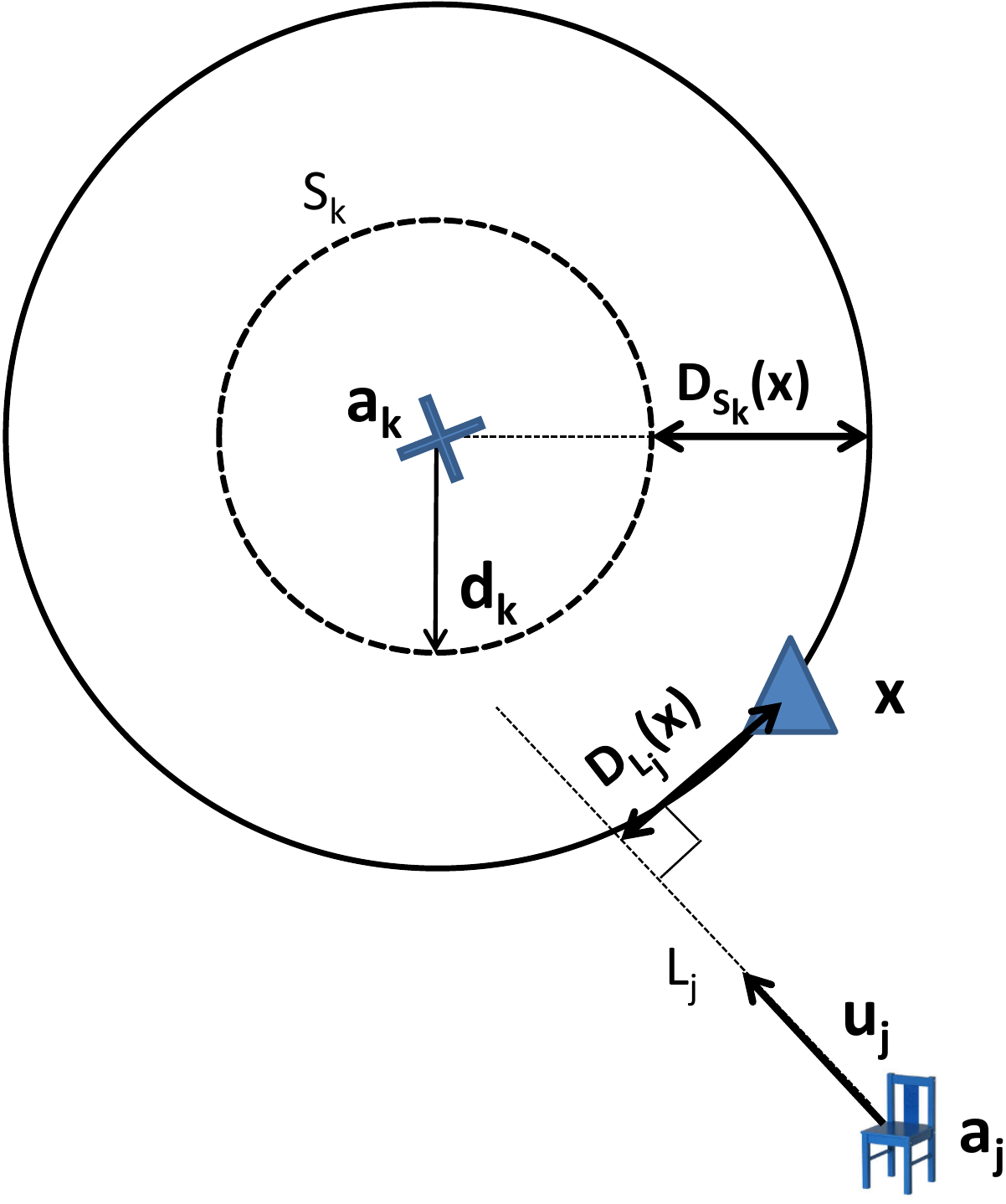}
  % \vspace{-4mm}
  \caption{Geometric representation of terms in the hybrid cost function \eqref{eq:FLORIS_cost_function}.}
  \label{fig:scheme_formulation}
  \vspace{-2mm}
\end{figure}

The ML formulation for range-only measurements corrupted by white Gaussian noise leads to a LS cost function with a sum of terms $(\| x - a_k\| - d_k)^2$. Each term may be interpreted as the squared distance of $x$ to a sphere of radius $d_k$ centered at $a_k$~\cite{SLCP:pinar:2014}, which penalizes deviations of $x$ from that sphere. Equivalently, this is the squared distance of $x - a_k$ to a sphere $S_k$ of radius $d_k$ centered at the origin, and will be denoted by $D_{S_k}^2(x -a_k)$ below. In the same spirit, our proposed approach for hybrid localization adds new terms to the cost function which penalize deviations of $x$ from lines with orientation $u_j$ starting at anchors $a_j$, analytically given by $(x - a_j)^T(\mathbf{I}_n - u_ju_j^T)(x - a_j)$~\cite{distancestolines:ballantine:1952}. Equivalently, this is the squared distance of $x - a_j$ to a line with orientation $u_j$ starting at the origin, and will be denoted by $D_{L_j}^2(x - a_j)$ (see Figure \ref{fig:scheme_formulation}). The cost function for single-source localization is 
\begin{equation}
  f(x) = \sum_{k \in \cal{R}} D_{S_k}^2(x - a_k) + \sum_{j \in \cal{T}} D_{L_j}^2(x - a_j).
  \label{eq:FLORIS_cost_function}
\end{equation}
The intuitive idea behind \eqref{eq:FLORIS_cost_function} is that this formulation attempts to balance, on the one hand, the distances of the target position estimate relative to spheres centered at selected wireless anchors with radii $d_k$ and, on the other hand, the distances to the lines originating at the other visual anchors with orientation $u_j$. Contrary to the range-only case, \eqref{eq:FLORIS_cost_function} is not equivalent to a likelihood function for common noise models in angular measurements $u_j$ (e.g., von Mises-Fisher). Despite its suboptimality, this LS formulation is taken as a starting point for the development of our algorithms due to its mathematical tractability.

In the cooperative scenario the heterogeneous WSN is represented by an undirected graph $\mathcal{G} = (\mathcal{V}, \mathcal{E})$, with the sensors as nodes and arcs representing the existence of pairwise measurements. For each unknown position of sensor node $i \in \mathcal{V}$, denoted by $x_i \in \mathbb{R}^n$, we define $\mathcal{A}_i$ as the subset of anchors whose distance and/or orientation to $i$ is available. Whenever necessary, $\mathcal{A}_i$ will be further split into $\mathcal{R}_i$ and $\mathcal{T}_i$, depending on which specific type of measurement --- range or angle --- pertains to each anchor.

We denote by $i \sim j$ the existence of measurements involving nodes $i$ and $j$, and by $N_{i}$ the set of neighbouring nodes of node $i$ in the graph. When necessary, these are more specifically denoted as range measurements, $i \stackrel{\mathcal{R}}{\sim} j$, or angle measurements $i \stackrel{\mathcal{T}}{\sim} j$, associated with restricted neighbourhoods $N^{\mathcal{R}}_{i}$ and $N^{\mathcal{T}}_{i}$. Actual measurements taken at node $i$ relative to node or anchor $j$ are denoted by $d_{ij}$ and $u_{ij}$, similarly to the single-source case. We adopt the following extension of the hybrid cost function \eqref{eq:FLORIS_cost_function} for cooperative localization
\begin{equation}
\label{eq:Cost_function_colab_fusion}
\begin{aligned}
%& \underset{\mathbf {x}}{\text{minimize}} 
& f(\mathbf{x}) = \sum_{i\stackrel{\mathcal{R}}{\sim} j} D_{S_{ij}}^2 \left(x_i - x_j \right) + \sum_{i} \sum_{k \in \mathcal{R}_i} D_{S_{ik}}^2 \left(x_i - a_{k} \right) \\
& + \sum_{i\stackrel{\mathcal{T}}{\sim} j} D_{L_{ij}}^2 \left(x_i - x_j \right) + \sum_{i} \sum_{k \in \mathcal{T}_{i}} D_{L_{ik}}^2 \left(x_i - a_{k} \right),
\end{aligned}
\end{equation}
where the argument $\mathbf{x}$ of \eqref{eq:Cost_function_colab_fusion} denotes the concatenation of all unknown sensor coordinates.
% In \eqref{eq:Cost_function_colab_fusion} $D_{S_{ij}}^2 \left(x_i - x_j \right)$ and $D_{S_{ik}}^2 \left(x_i - a_k\right)$ represent the squared euclidean distance operator (defined as $D_W^2(x) = \text{inf}_{y \in W} \|x - y \|^2$) of points $x_i - x_j$ and $x_i - a_k$ to the sets $S_{ij}$ and $S_{ik}$, spheres centered at the origin with radius $d_{ij}$ and $d_{ik}$, respectively. Similarly, $D_{L_{ij}}^2 \left(x_i - x_j \right)$ and $D_{L_{ik}}^2 \left(x_i - a_k\right)$ are the squared euclidean distances of points $x_i - x_j$ and $x_i - a_k$ to the straight lines $L_{ij}$ and $L_{ik}$ passing by the origin with orientation $u_{ij}$ and $u_{ik}$, respectively.  

In what follows, we will introduce algorithms that (approximately) solve both the source and network localization problems through convex relaxations of the hard original (nonconvex) functions \eqref{eq:FLORIS_cost_function}, \eqref{eq:Cost_function_colab_fusion}. The algorithms approximate the solution to estimate the sensor positions \{$x_i: i \in \mathcal{V}$\} ($\left|\mathcal{V}\right| = 1$ in the source localization case) from available noisy hybrid measurements. We follow different approaches for the single-source and collaborative cases, the former being cast as an SDP and solved in one shot, while the latter is tackled with a disk relaxation and a parallel and fast iterative algorithm that preserves the original optimization variables\footnote{Naturally, the collaborative algorithm can be used for single-source localization as well. This is illustrated in our simulation results of Section \ref{sec:results}.}.

\section{Hybrid Source Localization: FLORIS}
\label{sec:single-source}

Minimizing the single-source cost function \eqref{eq:FLORIS_cost_function} should find a source position $x$ that is close to the spheres $S_k$ centered at the wireless anchors, as well as the lines $L_j$ emanating from the visual anchors. Because the distance to any set $W$ from a point $x$ satisfies $D_W^2(x) = \inf_{y \in W} \| x - y\|^2$, one might consider a reformulation for the minimization of \eqref{eq:FLORIS_cost_function} where the ``closest points'' to $x$ in the sets $S_k$, $L_j$ are explicitly found. This yields the equivalent problem~\cite{Ferreira:hybrid_local:2015}
% \vspace{-3mm}
\begin{equation}
\label{eq:source_local_fusion}
\begin{aligned}
& \underset{x,y_k,\theta_k,t_k}{\text{minimize}}
& & \sum_{k=1}^{m} \| x - y_k \|^2  \\
& \text{subject to}
& & y_k = a_k + d_k \theta_k, \quad \| \theta_{k} \| = 1,  \quad k \in \cal{R}, \\
&&& y_k = a_k + u_k t_k, \quad t_k \in \mathbb{R}^+, \quad k \in \cal{T}.
% \vspace{-2mm}
\end{aligned}
\end{equation}
The first and second sets of constraints ensure that the ``closest points'' to the source, encoded by variables $y_k$, are located on the spheres $S_k$ or lines $L_k$, respectively. Then, \eqref{eq:source_local_fusion} jointly chooses $x$ and the $y_k$ to attain the best match (i.e., lowest dispersion). This is the same type of reformulation used for range-only source localization with the SLCP/SLNN algorithms~\cite{SLCP:pinar:2014}, but with new constraints for angular measurements.

Given all $y_k$, minimizing the cost function of \eqref{eq:source_local_fusion} is a standard least-squares problem whose optimal solution for $x$ is just the center of mass of the constellation $x = \frac{1}{m}\sum_{k} y_k$. This can be substituted back in the cost function to yield $\mathbf{y}^T\mathbf{J}\mathbf{y}$, where $\mathbf{y}$ is a vector of size $mn \times 1$ that stacks $y_1, \ldots, y_m$, and $\mathbf{J}$ is the projector onto the orthogonal complement of $\mathbf{1}_m\otimes\mathbf{I}_n = [\underbrace{\mathbf{I}_n \, \ldots \, \mathbf{I}_n}_m]^T$.
This can be compacted into matrix form as
%\vspace{-3ex}
\begin{equation}
\label{eq:joint_constraints}
\begin{aligned}
\mathbf{y} = \mathbf{a} + \mathbf{R} 
\begin{bmatrix} 
\boldsymbol{\theta} \\
\mathbf{t}\\
\end{bmatrix},
\end{aligned}
\end{equation}  
where $\boldsymbol{\theta}$ stacks the unit vectors $\theta_k \in \mathbb{R}^n$, $k \in \mathcal{R}$, and $\mathbf{t}$ stacks the scaling factors $t_k$, $k \in \mathcal{T}$. When $\mathcal{R} = \{ 1, \, \ldots, |\mathcal{R}| \}$ and $\mathcal{T} = \{ |\mathcal{R}|+1, \, \ldots, m \}$, the block diagonal-like matrix $\mathbf{R}$, of size $mn \times p$, with $p = n|\mathcal{R}| + |\mathcal{T}|$, is given by
\begin{equation}
  \label{eq:floris_R_mtx}
\mathbf{R} = \begin{bmatrix}
    d_1\mathbf{I}_n \\
    & \ddots \\
    & & d_{|\mathcal{R}|}\mathbf{I}_n \\
    & & & u_{|\mathcal{R}|+1} \\
    & & & & \ddots \\
    & & & & & u_m
  \end{bmatrix}.
\end{equation}
Problem \eqref{eq:source_local_fusion} is thus reformulated as
%\vspace{-1mm}
\begin{equation}
\label{eq:fusion_formulation}
\begin{aligned}
& \underset{\boldsymbol{\theta}, \mathbf{t}}{\text{minimize}}
& & \left( \mathbf{a} + \mathbf{R}
\begin{bmatrix} 
\boldsymbol{\theta} \\
\mathbf{t}\\
\end{bmatrix}
\right)^T \mathbf{J} \left( \mathbf{a} + \mathbf{R}
\begin{bmatrix} 
\boldsymbol{\theta} \\
\mathbf{t}\\
\end{bmatrix}
\right) \\
& \text{subject to}
& & \| \theta_k \| = 1, \; k \in \mathcal{R}, \qquad \mbox{$t_k \geq 0, \; k \in \mathcal{T}$}. 
\end{aligned}
\end{equation}
Cost function \eqref{eq:fusion_formulation} may be written as a quadratic form in $\begin{bmatrix}\boldsymbol{\theta}^T & \mathbf{t}^T & 1 \end{bmatrix}^{T}$, and expressed using the trace operator as
%\vspace{-2mm}
\begin{equation}
\label{eq:fusion_cost_trace}
\text{tr} \Bigl(
\underbrace{\begin{bmatrix} 
\mathbf{R}^T \mathbf{J} \mathbf{R} & \mathbf{R}^T\mathbf{J}\mathbf{a} \\
\mathbf{a}^T\mathbf{J} \mathbf{R} & \mathbf{a}^T \mathbf{J} \mathbf{a} \\
\end{bmatrix}}_{\mathbf{M}}
\underbrace{\begin{bmatrix}\boldsymbol{\theta}\\ \mathbf{t}\\ 1 \end{bmatrix}\begin{bmatrix}\boldsymbol{\theta}^T & \mathbf{t}^T & 1 \end{bmatrix}}_{\mathbf{W}} \Bigr).
\end{equation}
%\vspace{-1mm}
We now redefine matrix $\mathbf{W}$ as the optimization variable and rewrite the problem, yielding
\begin{equation}
\label{eq:fusion_formulation_matrix_final} 
\begin{aligned}
& \underset{\mathbf {W}}{\text{minimize}}
& & \text{tr}(\mathbf{M}\mathbf{W}) \\
& \text{subject to}
& & \mathbf{W} \succeq 0, \quad \text{rank}(\mathbf{W}) = 1\\
&&& \text{tr}(\mathbf{W}_{[k,k]}) = 1, \quad k \in \mathcal{R} \\
&&&  w_{[k],p+1} \geq 0, \quad k \in \mathcal{T} \\
&&& w_{p+1,p+1} = 1.
\end{aligned}
\end{equation}
The third constraint, where $\mathbf{W}_{[k,k]}$ denotes the submatrix of $\mathbf{W}$ comprising the rows/columns that pertain to $\theta_k$ in \eqref{eq:joint_constraints}, encodes $\| \theta_k \| = 1$. The fourth constraint, where $w_{[k],p+1}$ denotes the element of $\mathbf{W}$ in the row pertaining to $t_k$ and in the rightmost column (repeated in the bottom row of $\mathbf{W}$ due to symmetry), encodes $t_k \geq 0$. For $\mathcal{R}$ and $\mathcal{T}$ chosen as in \eqref{eq:floris_R_mtx} the pertinent elements of $\mathbf{W}$ are as follows
\begin{equation}
  \mathbf{W} = \begin{bmatrix}
    \mathbf{W}_{[1,1]} \\
    & \ddots \\
    & & \mathbf{W}_{[|\mathcal{R}|,|\mathcal{R}|]} & \ldots \\
    & & & & w_{[|\mathcal{R}|+1],p+1} \\
    & & & & \vdots \\
    & & & & w_{[m],p+1} \\
    & & & & w_{p+1,p+1}
  \end{bmatrix}.
\end{equation}
Finally, we drop the rank constraint in \eqref{eq:fusion_formulation_matrix_final} to obtain the relaxed SDP. Vectors $\boldsymbol{\theta}$ and $\mathbf{t}$ are obtained by singular value decomposition (SVD) factorization of the solution $\mathbf{W}$ or directly from its rightmost column (or bottom row), from which the $y_k$ are computed by \eqref{eq:joint_constraints} and the source position estimated as the average of these $m$ points.

As the complexity of SDPs in the class of interest is known to scale faster than cubically with the problem size~\cite{SLCP:pinar:2014}, it may be practically relevant to reduce the number of unknowns. Noting that the positivity constraints on the variables $t_k$ associated with the angles $u_k$ will usually be fulfilled naturally in all but extremely noisy situations (i.e., in our setup the source is very unlikely to be erroneously placed on ``the wrong side'' of visual markers), one may choose to drop those constraints. Then the quadratic problem may be solved in closed form for the unconstrained variables $t_k$ in terms of the $\theta_k$, eliminating them from the optimization problem. Appendix A provides the details of the reduced-complexity algorithm, whose accuracy is virtually identical to that of the formulation described above.

\section{Hybrid Sensor Network Localization: CLORIS}
\label{sec:colab}
%\subsection{Problem Formulation}
%\subsection{CLORIS}
Popular approaches for sensor network localization using range measurements cast the problem into the framework of Euclidean distance matrices, or related SDP formulations, and replace the unknown sensor coordinates with relaxed matrix variables. While (squared) distances between pairs of sensors may be elegantly expressed as linear constraints in terms of these relaxed variables, accounting for angular information is less natural. Here, we follow an alternative convex relaxation approach that preserves sensor coordinates and makes it easy to incorporate angular constraints.

Rather than using a ``one-shot'' relaxation for the hybrid cost function \eqref{eq:FLORIS_cost_function}, as in FLORIS, our approach here takes advantage of previous work which keeps the structure and variables of the extended (i.e., cooperative) cost function \eqref{eq:Cost_function_colab_fusion}. Moreover, this different approach can be solved with a lightweight and fast-converging gradient-like iterative algorithm that is easily parallelizable, unlike the ``one-shot'' SDP relaxation used in FLORIS. The proposed algorithm builds upon the disk relaxation developed in~\cite{soares:coopdistlocal:2014} for range-only network localization. The initial cost function in~\cite{soares:coopdistlocal:2014} resembles \eqref{eq:Cost_function_colab_fusion}, but only contains the ``circular'' squared distance terms $D_{S_{ij}}^2$, $D_{S_{ik}}^2$. These terms are nonconvex due to nonconvexity of the underlying sets $S_{ij}$, $S_{ik}$ (spheres), and the relaxation approach simply replaces the sets with their convex hulls (balls) $B_{ij}$, $B_{ik}$, yielding distances $D_{B_{ij}}^2 \left(z\right) = \text{inf}_{\| y \| \leq d_{ij}} \| z - y \|^2$ and similarly for source-anchor range measurements. The resulting cost function is \emph{not} the convex envelope of the original nonconvex cost\footnote{Even though $D_{B_{ij}}^2$, $D_{B_{ik}}^2$ \emph{are} the convex envelopes of $D_{S_{ij}}^2$, $D_{S_{ik}}^2$~\cite{soares:coopdistlocal:2014}.}, but as an underestimator it does provide remarkably good approximation even with few sensors/anchors, with quantifiable sub-optimality. 

Squared distances to lines $D_{L_{ij}}^2$, $D_{L_{ik}}^2$ in \eqref{eq:Cost_function_colab_fusion} are already convex, and can be directly inserted into the relaxed hybrid cost function. CLORIS, our algorithm for cooperative localization, thus solves the convex problem
\begin{equation}
\label{eq:Cost_function_colab_fusion_final}
\underset{\mathbf {x}}{\text{minimize}} \hat{f}(\mathbf{x}) = g(\mathbf{x}) + h(\mathbf{x})
\end{equation}
\begin{gather*}
g(\mathbf{x}) = \sum_{i\stackrel{\mathcal{R}}{\sim} j} {\textstyle \frac{1}{2}}D_{B_{ij}}^2 \left(x_i - x_j \right)  + \sum_{i\stackrel{\mathcal{T}}{\sim} j} {\textstyle \frac{1}{2}}D_{L_{ij}}^2 \left(x_i - x_j \right)\\
h(\mathbf{x}) = \sum_{i} \left\{ \sum_{k \in \mathcal{R}_i} {\textstyle \frac{1}{2}}D_{B_{ik}}^2 \left(x_i - a_{k} \right) + \sum_{k \in \mathcal{T}_{i}} {\textstyle \frac{1}{2}}D_{L_{ik}}^2 \left(x_i - a_{k} \right) \right\}.
\end{gather*}
The factor ${\textstyle \frac{1}{2}}$ is prepended to all squared distance terms above to streamline the expressions for the gradients. The strategy to efficiently obtain the solution of \eqref{eq:Cost_function_colab_fusion_final} parallels the cooperative synchronous algorithm of~\cite{soares:coopdistlocal:2014} based on Nesterov's accelerated gradient descent, introduced in~\cite{Nesterov:gradient:1983}. The computational simplicity of each iteration stems from the focus on first-order methods, while fast convergence follows directly from known optimality properties of Nesterov's technique for this class of algorithms, which assert that no other descent gradient method can perform better. 
Recently,~\cite{Piovesan:convex-nonconvex:2016} has used the convex relaxation of~\cite{soares:coopdistlocal:2014} as a first stage and then transitions to the nonconvex non-relaxed formulation to force the solution to adhere to the original nonconvex problem. The localization accuracy of this approach~\cite{Piovesan:convex-nonconvex:2016} is equivalent to the one in~\cite{soares:coopdistlocal:2014}, however~\cite{Piovesan:convex-nonconvex:2016} achieves a faster convergence speed, due to the use of ADMM instead of Nesterov method. 
Nonetheless, we focus our attention on algorithms that can be solved in a distributed manner, and resorting to an ADMM method introduces a communication overhead, as its primal-dual approach creates several auxiliary variables that must be passed across nodes. Hence, despite the advance of~\cite{Piovesan:convex-nonconvex:2016}, we opt to follow the original approach from~\cite{soares:coopdistlocal:2014} that reduces the number of communications in a distributed scenario, even at the cost of slower convergence.    
Simulation results show that the approach of~\cite{soares:coopdistlocal:2014} matches and often exceeds the accuracy of state-of-the-art range-only network localization algorithms at a fraction of their computational (and communication, if applicable) complexity.

Nesterov's method is applicable since the cost function $\hat{f}(\mathbf{x})$ is differentiable and both terms that comprise the gradient, $\nabla g(\mathbf{x})$ and $\nabla h(\mathbf{x})$, are Lipschitz continuous, with an easily computable Lipschitz constant $L_{\hat{f}} = L_{g} + L_{h}$. Although all sets in the cost function of~\cite{soares:coopdistlocal:2014} are balls, the expressions for the Lipschitz constant given there are still valid when squared distances pertain to other convex sets, including the lines $L_{ij}$, $L_{ik}$ in \eqref{eq:Cost_function_colab_fusion_final}, and require no modification\footnote{The upper bound for the Lipschitz constant derived in~\cite{soares:coopdistlocal:2014} depends on the maximum degree of a node (number of neighbouring nodes) and maximum number of visible anchors. The result does not depend on the specifics of the projection operator, and thus remains valid here when the full graph of pairwise range or angular measurements and full set of anchors is considered.}.

Importantly, the algorithms proposed in~\cite{soares:coopdistlocal:2014} take advantage of the special form of the gradient of the relaxed cost function to compute it in a distributed manner, through local computations at each node and synchronous or asynchronous exchanges with its neighbours. This relevant property from~\cite{soares:coopdistlocal:2014} does carry over to CLORIS since, as discussed below, the gradients of the new terms $D_{L_{ij}}^2$, $D_{L_{ik}}^2$ can still be evaluated in a distributed manner, similarly to those involving range measurements. In analogy with~\cite{soares:coopdistlocal:2014} CLORIS is thus described below in terms of local computation and communication between nodes, even though the parallel algorithm may actually end up being executed at a central processor, either sequentially or mapping ``nodes'' as concurrent threads in a parallel implementation.

In addition to algorithmic aspects, a truly distributed solution might reasonably be expected to meet the requirement that measurements for node $i$ or pair $i \sim j$ be obtained locally. This is easily satisfied for ranges and node-anchor angles, but not for internode angular measurements unless node orientations (attitudes) are known using, e.g., on-board inertial measurement units (IMU). Here we emphasize scenarios where the assumed hardware simplicity of nodes excludes such devices, and incidence directions $u_{ij}$ are derived from pairwise translation and rotation matrices (see Section \ref{sec:self-calib}), recursively propagated across the node visibility subgraph starting at anchors. Therefore, depending on the measurement hardware, setting up pairwise measurements for CLORIS might require a level of coordination across the network that is not apparent in the algorithmic description.

At a high level, each node position estimate in CLORIS is initialized arbitrarily and its current value at time $k$, $x_i(k)$, is updated resorting to a combination $w_i$ of the two previous iterations ($x_i(k-1)$ and $x_i(k-2)$). The processing loop can be performed up to a predefined maximum number of iterations, or it can stop when the gradient is smaller than a preset threshold (the latter is computed centrally, and therefore does not qualify as a native distributed stopping criterion). The pseudo-code of CLORIS is presented in Algorithm~\ref{alg:Parallel_method}. 
\begin{algorithm}
 \caption{CLORIS: Parallel hybrid network localization based on accelerated Nesterov's optimal method 
   \label{alg:Parallel_method}}
  \begin{algorithmic}[1]
\Require $L_{\hat{f}}$, \{$d_{ij}, \, u_{ij} : i \sim j \in \mathcal{E}$\}, \{$d_{ik}, \, u_{ik} : i \in \mathcal{V}, k \in \mathcal{A}$\}; 
\Ensure $\mathbf{\hat{x}}$;
\item $k=0$;
\item $x_i(0)$ and $x_i(-1)$ is randomly assigned, for each node $i$;
\While{stopping criterion is not met, each node $i$}
	\State $k=k+1$;
	\State $w_i = x_i(k-1) + \frac{k-2}{k+1}\left(x_i(k-1) - x_i(k-2)\right)$;
	\State $w_i$ is broadcast to the neighbours of node $i$;
	\State \label{eq:gradg} \vspace{-3ex}\begin{equation*}\begin{split} & \nabla g(w_i) = \\ & |N^{\mathcal{R}}_i| w_i - \sum\limits_{j \in N^{\mathcal{R}}_{i}} w_j + \sum\limits_{j \in N^{\mathcal{R}}_{i}} c^{\mathcal{R}}_{ji}{P_{B_{ij}}} \left(w_i - w_j\right) + \\ & |N^{\mathcal{T}}_i| w_i - \sum\limits_{j \in N^{\mathcal{T}}_{i}} w_j + \sum\limits_{j \in N^{\mathcal{T}}_{i}} c^{\mathcal{T}}_{ji}{P_{L_{ij}}} \left(w_i - w_j\right); \end{split}\end{equation*}
	\State \label{eq:gradh} \vspace{-3ex}\begin{equation*}\begin{split} & \nabla h(w_i) = \\ & |\mathcal{A}_i|w_i - \sum\limits_{k \in \mathcal{R}_i} P_{B_{ik}}(w_i - a_k) - \sum\limits_{k \in \mathcal{T}_i} P_{L_{ik}}(w_i - a_k); \end{split}\end{equation*}
	\State $x_i(k) = w_i - \frac{1}{L_{\hat{f}}} \left(\nabla g_i (w_i) + \nabla h_i(w_i)\right)$;
\EndWhile 
\State \textbf{return} $\mathbf{\hat{x}} = \mathbf{x}(k)$
  \end{algorithmic}
\end{algorithm}
In Algorithm~\ref{alg:Parallel_method}, $P_{B_{ij}}(z)$ denotes the orthogonal projection of point $z$ onto the ball $B_{ij}$, while $P_{L_{ij}}(z) = u_{ij}u_{ij}^Tz$ is the orthogonal projection of $z$ onto the straight line $L_{ij}$. The weight $c^{\mathcal{R}}_{ji}$ denotes an element of any directed arc-node incidence matrix (agreed upon at the onset of the algorithm) for the subgraph of internode range measurements that pertains to one of the edges connected to node $i$; similarly for $c^{\mathcal{T}}_{ji}$ and the subgraph of internode angle measurements.

Focusing on the new angle-related terms that were absent in the derivation of the distributed algorithm in~\cite{soares:coopdistlocal:2014}, we first recall that the gradient of the squared distance to a convex set $W$ satisfies
\begin{equation}
\nabla  {\textstyle \frac{1}{2}}D_{W}^2(z) = z - P_W(z),
\end{equation}
where $P_W(z)$ is the orthogonal projection of $z$ onto $W$. For node-anchor angular measurements the contribution to the components of $\nabla h(\mathbf{x})$ in \eqref{eq:Cost_function_colab_fusion_final} related to $x_i$ is thus
\begin{equation}
\sum_{k \in \mathcal{T}_{i}} x_i - P_{L_{ik}}(x_i - a_k) = |\mathcal{T}_{i}|x_i - \sum_{k \in \mathcal{T}_{i}} P_{L_{ik}}(x_i - a_k),
\end{equation}
with $P_{L_{ik}}(x_i - a_k) = u_{ik}u_{ik}^T(x_i - a_k)$. When added to the contributions from range measurements, this yields line \ref{eq:gradh} in Algorithm~\ref{alg:Parallel_method}.

For internode measurements let us define $D_L^2$ as the squared distance to the Cartesian product of all lines $L_{ij}$, indexed by the concatenation of all pairwise differences $x_i - x_j$, $i\stackrel{\mathcal{T}}{\sim} j$. We write the latter compactly as $\mathbf{A}\mathbf{x}$ in terms of the vector of concatenated sensor coordinates $\mathbf{x}$ and matrix $\mathbf{A}$, obtained from a directed arc-node incidence matrix $\mathbf{C}^{\mathcal{T}}$ for the subgraph of angle measurements as $\mathbf{A} = \mathbf{C}^{\mathcal{T}} \otimes \mathbf{I}_p$. Then, the contribution of these terms to $g(\mathbf{x})$ is ${\textstyle \frac{1}{2}}D_{L}^2 \left( \mathbf{A}\mathbf{x} \right) \stackrel{\Delta}{=} \sum_{i\stackrel{\mathcal{T}}{\sim} j} {\textstyle \frac{1}{2}}D_{L_{ij}}^2 \left(x_i - x_j \right)$, with gradient~\cite{Magnus2007}
\begin{equation}
  \begin{split}
    \nabla {\textstyle \frac{1}{2}}D_{L}^2 \left( \mathbf{A}\mathbf{x} \right) & = \mathbf{A}^T\left( \mathbf{A}\mathbf{x} - P_L\left( \mathbf{A}\mathbf{x} \right) \right) \\ & = \bm{\mathcal{L}}\mathbf{x} - \mathbf{A}^T P_L\left( \mathbf{A}\mathbf{x} \right),
  \end{split}
\end{equation}
where $\bm{\mathcal{L}}$ is the Laplacian matrix of the subgraph and $P_L$ denotes the concatenation of projections onto the lines $L_{ij}$. This expression, written explicitly in terms of individual sensor coordinates, is added to range-related terms in line \ref{eq:gradg} of Algorithm~\ref{alg:Parallel_method}.

\section{Self-calibration for assimilation of range and video sensory data}
\label{sec:self-calib}

Recall from Figure \ref{fig:scheme_formulation} that range and bearing are captured by different, not necessarily collocated, sensors. To be used in our hybrid cost functions based on Euclidean distances, the spatial information gained from these measurements has to be expressed in a common reference frame, which we choose as any convenient one for describing the coordinates of anchors in a particular setup. Angular information can be provided by a number of technologies, but the emphasis here is on low-cost solutions based on video and fiducial markers such as ARUCO~\cite{Aruco2014}, where a supplied library function is invoked as a black box to robustly express the camera pose with respect to a fiducial marker\footnote{When multiple ARUCO markers are deployed, a documented calibration and preprocessing procedure may be carried out such that subsequent library calls return pose information relative to a fixed (reference) marker, regardless of which particular marker is actually detected in a given camera image. This operating mode is assumed in our derivations.}. It is then necessary to translate this to the global reference frame, which would be straightforward if the position and pose of the marker, as well as the relative positions of the range sensor and the camera, were exactly known. Manually determining these parameters is an option, but the process is cumbersome and error prone; e.g., range sensors are not point-like, and distances to the camera should be measured with respect to its focal point, which is usually not accessible as it is located inside the enclosure. A novel and more convenient alternative, described in this section, is to automate the calibration process through a ``system identification'' procedure.

\begin{figure}[h]
    \centering
    \includegraphics[width=0.5\columnwidth]{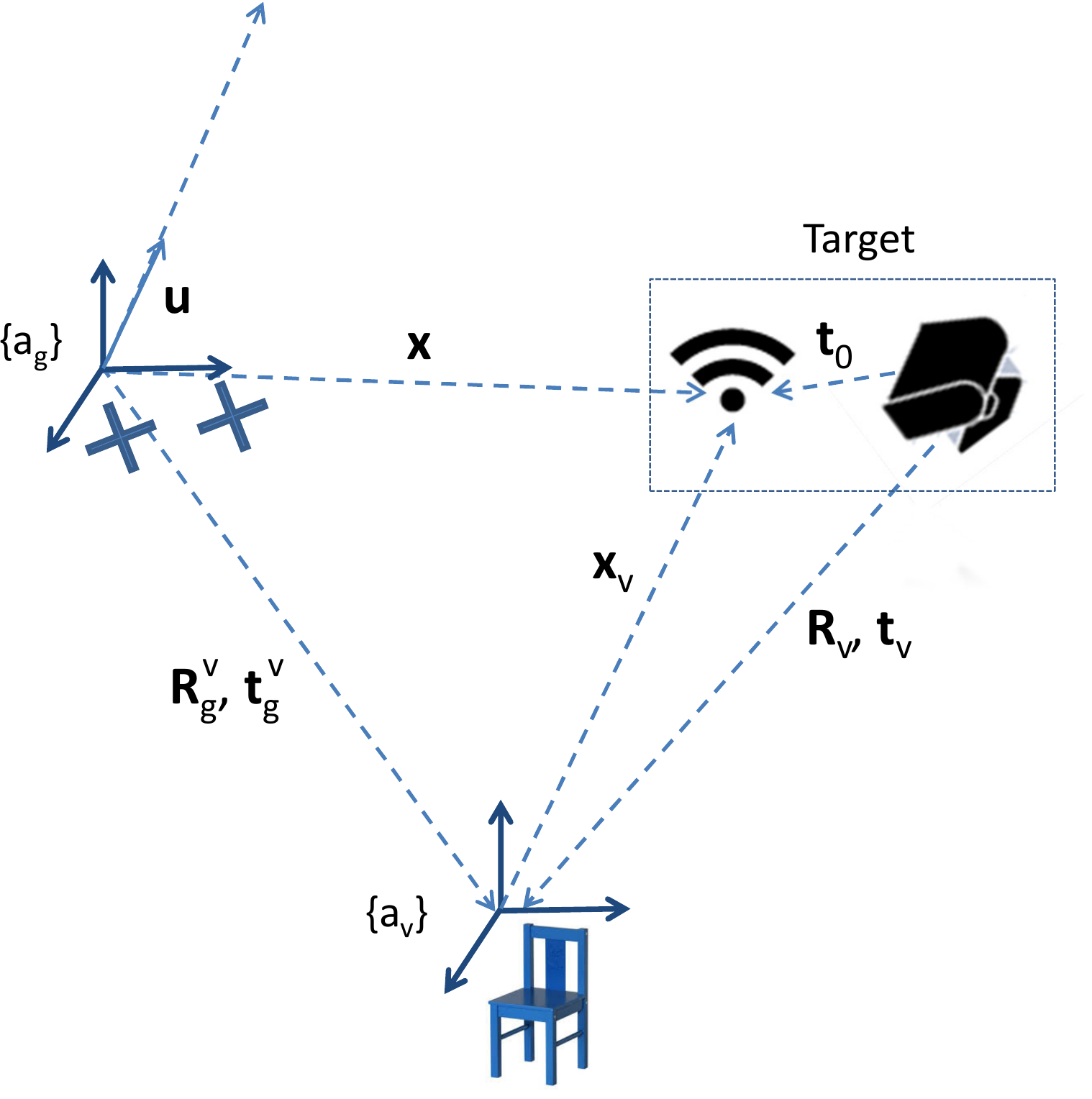}
    \caption{Scheme for the self-calibration of the two uncoupled networks, auxiliary to the derivation of expression \eqref{eq:auto-calib_syseq}. See also Figure \ref{fig:experimental_setup_glass} with the respective experimental set-up.}
    \label{fig:scheme_calib}
            % \vspace{-4mm}
  \end{figure}

Figure \ref{fig:scheme_calib} depicts the setup for a single target comprising a range sensor (e.g., acoustic or ultrawideband) at position $x$ in the global reference frame $\{a_g\}$ and a camera whose pose (rotation matrix $\mathbf{R}_\text{v}$ and translation vector $\mathbf{t}_\text{v}$) is measured relative to the local reference frame of a visual anchor $\{a_v\}$ (abbreviated as visual frame). To enable calibration, the target first evolves in regions where it can simultaneously determine its position using \emph{only} range measurements, and measure its pose relative to some of the visual anchors. The missing parameters are determined upon collecting a sufficiently large (say, 20) and diverse set of range-based location estimates, $x_{r}$, and visual poses, $(\mathbf{R}_\text{v}, \, \mathbf{t}_\text{v})$.
% \vspace{-4mm} 

The strategy for calibration is to express the position of the target's range sensor in the global frame through the composition of two rigid transformations involving measured poses, and matching it to the same position computed from range measurements only. Referring to Figure \ref{fig:scheme_calib}, the first rigid transformation, $(\mathbf{R}_g^v, \, \mathbf{t}_g^v)$, translates points in the visual frame into points in the global frame. The second one translates points from the visual frame into coordinates in the camera viewpoint. The unknowns are $\mathbf{R}_g^v$, $\mathbf{t}_g^v$, and translation vector $\mathbf{t}_\text{0}$, which reflects the position of the range sensor in the natural camera-centric 3D coordinate system.

The target position is represented in both frames by the rigid body equations
\begin{equation}
\label{eq:auto-calib_syseq}
\begin{split}
x_\text{v} & = \mathbf{R}_\text{v} \mathbf{t}_\text{0} + \mathbf{t}_\text{v},\\
 x & = \mathbf{R}_g^v  x_\text{v} + \mathbf{t}_g^v.
\end{split}
\end{equation}
The first expression sets the position of the range sensor in the visual frame by adding the translation $\mathbf{t}_\text{0}$, with appropriate rotation, to the position of the camera $\mathbf{t}_\text{v}$. The second expression maps $x_\text{v}$ to the global frame through a rigid body transformation. The bearing $u$ to be used in the hybrid cost function is a unit-norm vector with the same direction as $\mathbf{R}_g^vx_\text{v}$, i.e., $x_\text{v}$ rotated to the global frame.

We identify $\mathbf{R}_g^v$, $\mathbf{t}_g^v$, $\mathbf{t}_\text{0}$ from source (target) positions $x_{r}$ (computed via a range-based algorithm such as SR-LS~\cite{sourceLS:beck:2008} or SLNN~\cite{SLCP:pinar:2014}) and camera poses $(\mathbf{R}_\text{v}, \, \mathbf{t}_\text{v})$ through a nonlinear constrained least-squares procedure. In general, problems of this sort are very hard to solve numerically and prone to local minima. However, we show that we may decouple the estimation of $\mathbf{R}_g^v$, $\mathbf{t}_\text{0}$ from $\mathbf{t}_g^v$ by taking differences of range-based estimated positions from pairwise observations, leading to sub-problems that are solvable in closed form.
% These approximations are tight and convergence to the global minimum was always attained.

From \eqref{eq:auto-calib_syseq}, the following relation should ideally hold for the pairwise difference between estimated/observed quantities for measurements $i$ and $j$
\begin{equation}
\label{eq:auto-calib_tdiff_simp}
\underbrace{x_{r_i} - x_{r_j}}_{\tilde{x}_{r_{ij}}}  = \mathbf{R}_g^v \Bigl[ \; \underbrace{(\mathbf{R}_{\text{v}_i} - \mathbf{R}_{\text{v}_{j}})}_{\widetilde{\mathbf{R}}_{\text{v}_{ij}}}  \mathbf{t}_\text{0} + \underbrace{(\mathbf{t}_{\text{v}_i} - \mathbf{t}_{\text{v}_j})}_{\tilde{\mathbf{t}}_{\text{v}_{ij}}} \; \Bigr],% _\text{$\mathtbf{t_{cam_k}}$} \left],
\end{equation}
which does not depend on $\mathbf{t}_g^v$. In the presence of noise and localization errors the rotation $\mathbf{R}_g^v$ and translation $\mathbf{t}_\text{0}$ are estimated from multiple pairwise measurements by solving the LS problem
\begin{equation}
\label{eq:auto-calib_error_function}
\begin{aligned}
& \underset{\mathbf{R}_g^v, \mathbf{t}_\text{0}}{\text{minimize}} 
& & \sum_{i,j} \| \tilde{x}_{r_{ij}} - \mathbf{R}_g^v ( \widetilde{\mathbf{R}}_{\text{v}_{ij}} \mathbf{t}_\text{0} +  \tilde{\mathbf{t}}_{\text{v}_{ij}} ) \|^2 \\
& \text{subject to}
& & {\mathbf{R}_g^v}^T \mathbf{R}_g^v = \mathbf{I}_3. 
\end{aligned}
\end{equation}
Relative to each of the unknowns, \eqref{eq:auto-calib_error_function} can be independently minimized in closed form: if $\mathbf{t}_\text{0}$ is known, $\mathbf{R}_g^v$ is the solution of a Procrustes problem~\cite{matrixcomputs:Golub:1991}; and knowing $\mathbf{R}_g^v$, then $\mathbf{t}_\text{0}$ is given by a pseudoinverse. Inserting $\mathbf{R}_g^v$ and $\mathbf{t}_\text{0}$ back into~\eqref{eq:auto-calib_syseq}, $\mathbf{t}_g^v$ is computed in closed form as the average value of $x_{r_{i}} - \mathbf{R}_g^v  (\mathbf{R}_{\text{v}_i} \mathbf{t}_\text{0} + \mathbf{t}_{\text{v}_i})$ across measurements.

Algorithm \ref{alg:Self-calib} summarizes the steps of the calibration procedure, which is described and analyzed in detail in~\cite{Calibration_Report:Ferreira:2015}. 
%In the interest of space we omit the proof, given in [\textcolor{red}{ref technical report}], that the algorithm converges to a fixed point. 
In~\cite{Calibration_Report:Ferreira:2015} we also show, in simulation, that the true solution was consistently obtained with this algorithm. With real data, the solution was consistent across trials and compatible with the characteristics of our testbed.

\begin{algorithm}
 \caption{Self-calibration procedure to support conversion of visual sensor measurements to global angular data
   \label{alg:Self-calib}}
  \begin{algorithmic}[1]
\Require \{$x_{r_i}, \, \mathbf{R}_{\text{v}_i}, \, \mathbf{t}_{\text{v}_i}: i = 1, \ldots N_\text{calibration}$\}; 
\Ensure $\mathbf{R}_g^v$, $\mathbf{t}_g^v$, $\mathbf{t}_\text{0}$;
\item $k=0$;
\item Initialize $\mathbf{t}_\text{0} = 0$, $\varepsilon_{\mathbf{R}_g^v} = \infty$;
\item Choose set of pairwise differences \{$\tilde{x}_{r_{ij}}, \, \widetilde{\mathbf{R}}_{\text{v}_{ij}}, \, \tilde{\mathbf{t}}_{\text{v}_{ij}}$\};
\While{stopping criterion on $\varepsilon_{\mathbf{R}_g^v}$ is not met}
	\State $k=k+1$;
	\State Find $\mathbf{R}_g^v(k)$ by solving a Procrustes problem with \{$\tilde{x}_{r_{ij}}$\} and \{$\widetilde{\mathbf{R}}_{\text{v}_{ij}} \mathbf{t}_\text{0} + \tilde{\mathbf{t}}_{\text{v}_{ij}}$\} as the two sets of points;
	\State \vspace{-3.5ex}\begin{equation*}\hspace{3ex}\mathbf{t}_\text{0}(k) = \underset{\mathbf{t}_\text{0}}{\text{arg min}} \sum_{i,j} \| \tilde{x}_{r_{ij}} - \mathbf{R}_g^v(k) ( \widetilde{\mathbf{R}}_{\text{v}_{ij}} \mathbf{t}_\text{0} +  \tilde{\mathbf{t}}_{\text{v}_{ij}} ) \|^2;\end{equation*}
% \left(\sum_{i,j} \mathbf{R}_\textsubscript{k}}^T \mathbf{R}_\text{diff\textsubscript{k}} \right)^{-1} ( \sum_{k} \mathbf{R}_\text{diff\textsubscript{k}}^T \mathbf{R}_\text{net}(i)^T x_\text{diff\textsubscript{k}}$ \item [] $- \sum_{k} \mathbf{R}_\text{diff\textsubscript{k}}^T  \mathbf{T}_\text{v\textsubscript{k}})$;
	\State $\varepsilon_{\mathbf{R}_g^v} = \|\mathbf{R}_g^v(k) - \mathbf{R}_g^v(k-1)\|_{F}$;
\EndWhile 
\State $\mathbf{t}_g^v = \overline{x}_r - \mathbf{R}_g^v(k) \left(\overline{\mathbf{R}}_\text{v} \mathbf{t}_\text{0}(k) + \overline{\mathbf{t}}_\text{v} \right)$;
\State \textbf{return} $\mathbf{R}_g^v = \mathbf{R}_g^v(k)$, $\mathbf{t}_\text{0} = \mathbf{t}_\text{0}(k)$, $\mathbf{t}_g^v$
  \end{algorithmic}
\end{algorithm}

\section{Numerical and Experimental Results}
\label{sec:results}
  \begin{figure}[h]
    \centering
    \includegraphics[width=0.75\columnwidth]{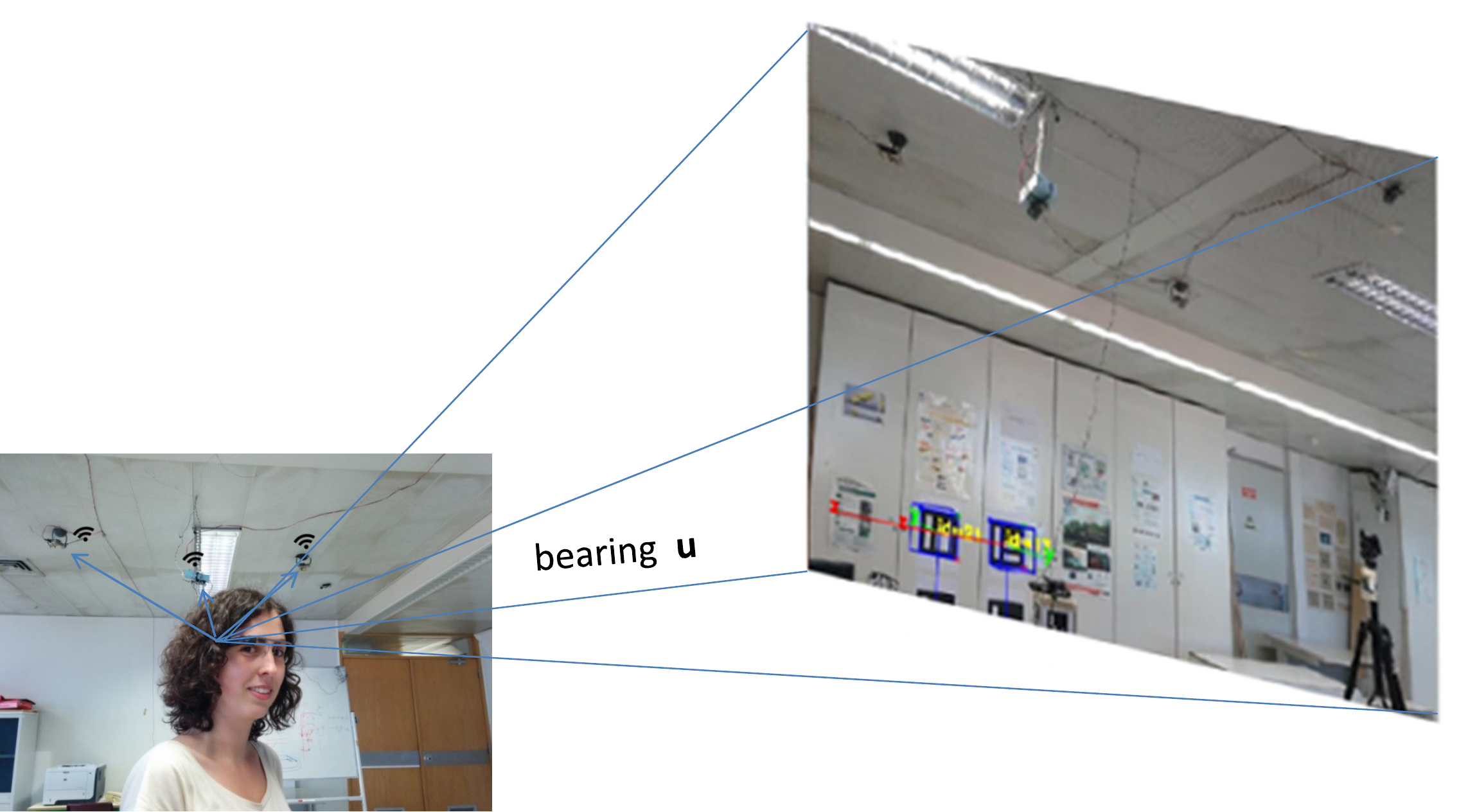}
    \caption{Experimental set-up for localization, with Cricket modules as ranging sensors and video camera from Google Glass detecting ARUCO tags as visual sensor.}
    \label{fig:experimental_setup_glass}
            \vspace{-2mm}
  \end{figure} 

In this section we evaluate the performance of the two newly proposed localization methods in simulation and experimentally.

Firstly, FLORIS is benchmarked against SR-LS~\cite{sourceLS:beck:2008} and SLNN~\cite{SLCP:pinar:2014} in simulation. The latter was chosen based on the assessment of~\cite{SLCP:pinar:2014} which showed that SLNN has higher accuracy in 3D than previously proposed SDP-based methods. While SR-LS does not belong to the SDP class, being somewhat less accurate due to the use of \emph{squared} ranges that amplify the effect of noise, it does provide an interesting trade-off between good precision and low computational complexity. Then, findings from implementing FLORIS in the deployed set-up shown in Figure \ref{fig:experimental_setup_glass} are also presented and discussed.

Our CLORIS distributed sensor network localization method is also simulated and tested experimentally. The attained results are compared, for both cases, against the range-only precursor algorithm based on a disk relaxation~\cite{soares:coopdistlocal:2014}, which was proven to be more efficient and accurate than previous range-based methods. Also, we benchmark CLORIS against the hybrid range/bearing SDP relaxation of~\cite{Biswas05integrationof}. Since this method is specific to 2D, we only compare it to CLORIS through simulation. 

\paragraph*{\emph{Methods}}
We test single-source localization algorithms in simulation on randomly-generated networks of acoustic and visual anchors in a unit square (in 2D) or a unit cube (in 3D). Pairwise measurement graphs that arise in collaborative scenarios are more prone to lead to ill-posed localization problems where a unique set of sensor positions consistent with the data cannot be determined~\cite{SensorNetTheoryGraph:anderson:2010}. Hence, we test CLORIS in simulation on a fixed set of pre-generated localizable networks, rather than generating random networks, testing for rigidity~\cite{Aspnes:TheoryNetLoc:2006} (difficult in 3D), and discarding non-localizable instances.

To simulate range and/or bearing measurements, we first add white Gaussian noise to differences of true position vectors according to $\delta = \delta_{0} + w$, where $\delta_{0}$ denotes a nominal sensor/sensor or sensor/anchor difference, viz.\ $x_{i} - x_{j}$ or $x_{i} - a_{k}$. The Gaussian noise term conforms to $w ~ \sim \mathcal{N}(0,\eta^2\| \delta_{0} \|^2\mathbf{I})$, where $\eta$ is referred to below as \emph{noise factor}. Then, we generate range and bearing measurements as $d = \| \delta \|$ and $u = \delta/\| \delta \|$. Under this model range errors tend to increase for longer distances, reflecting a behaviour often found in real-world range measurement systems.

Estimation accuracy is assessed by the Root-Mean-Square Error (RMSE). For a set of $MC$ Monte Carlo runs on a scenario with $N$ unknown node positions, it is defined as
\begin{equation}
  \label{eq:rmse}
  RMSE = \sqrt{\frac{1}{MC} \frac{1}{N} \sum_{k=1}^{MC} \sum_{i=1}^{N} \|x_i - \hat{x_i}^k\|^2 },
\end{equation}
where $x_i$ and $\hat{x_i}^k$ denote the true and estimated positions of the $i$-th node in the $k$-th Monte Carlo run, respectively. Sets of $MC=1000$ trials were run for each measurement noise factor $\eta \in (0.001, 0.005, 0.01, 0.05, 0.1, 0.2, 0.3, 0.4)$. Additionally, for the single-source paradigm, we evaluate the performance of FLORIS also based on the rank of matrix $\mathbf{W}$ in the relaxed reformulation of \eqref{eq:fusion_formulation_matrix_final}, which ideally should be one.
% When $\text{rank}(\mathbf{W}) = 1$, the relaxed solution coincides with that of the non-relaxed problem.

%\textcolor{red}{Another performance measure commonly used in estimation is the Cram\'{e}r-Rao Lower Bound (CRLB), which sets a lower bound on the variance of an unbiased estimator. Hence, we also show, for the network localization, the CRLB computed as $\sqrt{\frac{1}{N} \text{tr}(\mathbf{F}^{-1})}$, where $\mathbf{F}$ is the Fisher Information Matrix.}

The following experiments were run in a machine powered by an Intel Core i7-2600 CPU @ 3.40GHz and 8GB of RAM, using \emph{MATLAB R2013a} and the general-purpose SDP solver CVX/SDPT3, for solving \eqref{eq:fusion_formulation_matrix_final} in FLORIS and the SDP-based method of~\cite{Biswas05integrationof}.  

\subsection{Source Localization (FLORIS)}
\subsubsection{Numerical Results}
\begin{figure}[!h]
  \centering
  \subfloat[]{%
    \label{RMSE_floris_sim_comparison_2d_convert}%
    \includegraphics[width=0.5\linewidth]{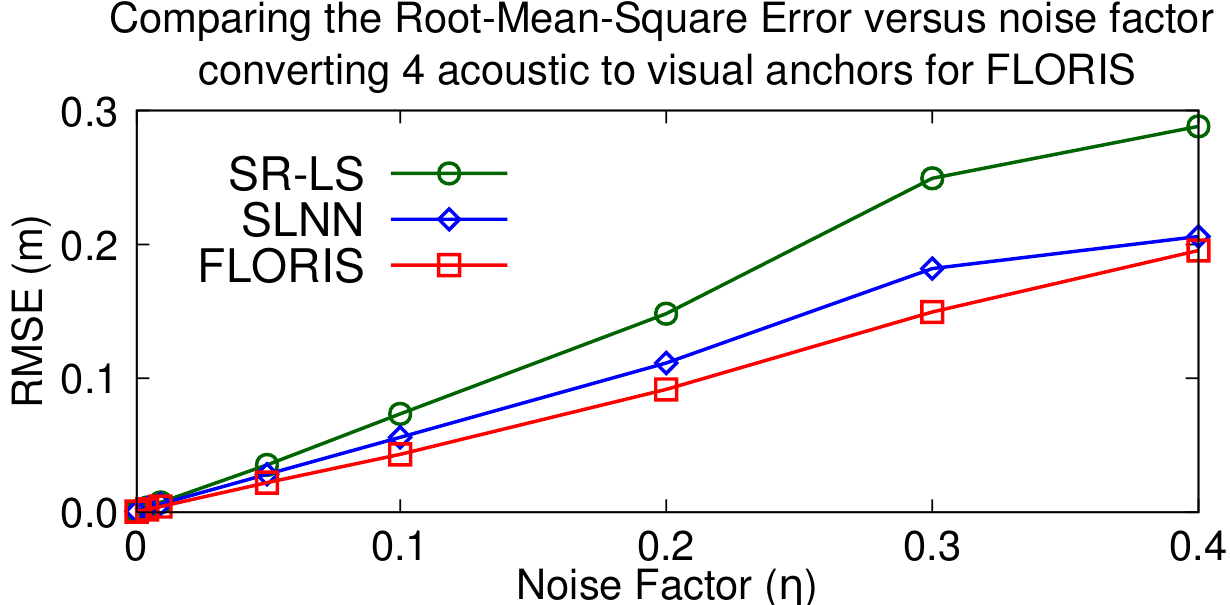}}
  \subfloat[]{%
    \label{RMSE_floris_sim_comparison_3d_convert}%
    \includegraphics[width=0.5\linewidth]{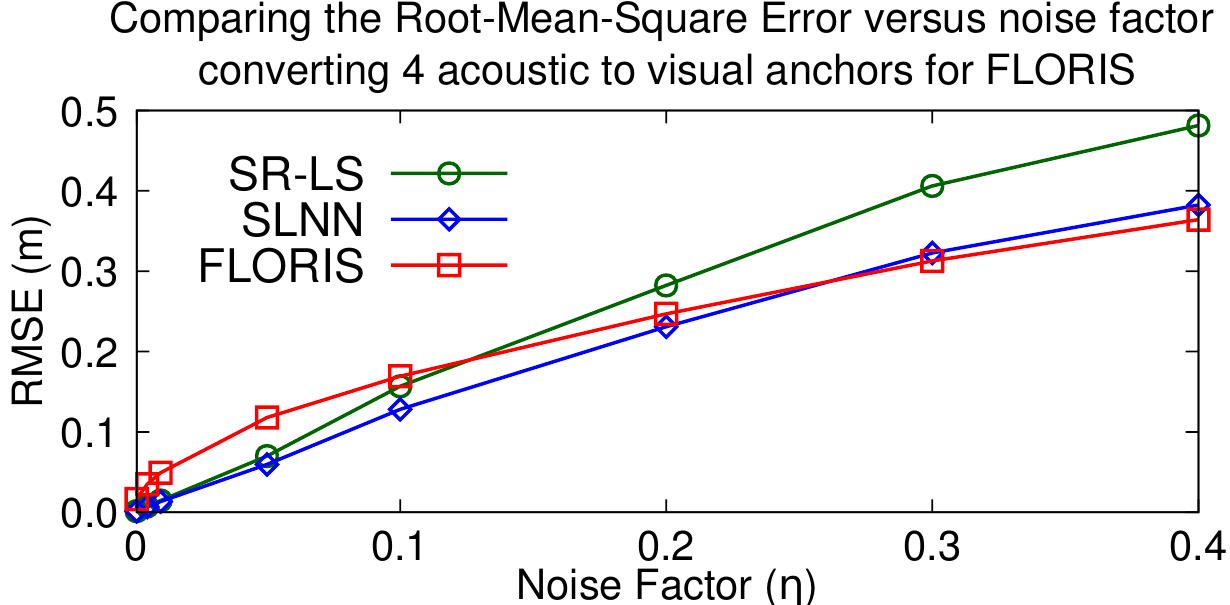}}
  \caption[]{Single-source localization: RMSE vs.\ noise factor for 8 acoustic or 4 acoustic + 4 visual anchors. \subref{RMSE_floris_sim_comparison_2d_convert} 2D. \subref{RMSE_floris_sim_comparison_3d_convert} 3D.}
  \label{RMSE_floris_sim_comparison_convert}
\end{figure}

~\\
\emph{Example 1:} Figure \ref{RMSE_floris_sim_comparison_convert} shows a performance comparison of the three methods in 2D and 3D for several values of noise factor $\eta$. The random network configurations comprise 8 acoustic anchors for SLNN and SR-LS. For FLORIS, 4 of these anchors were converted to visual ones, so that all algorithms are fairly compared based on the same total number of anchors. As expected, SR-LS does worse than the other algorithms due to squaring of noisy range measurements, particularly for higher noise factors. FLORIS consistently outperforms SLNN in 2D, albeit not by a large margin, whereas in 3D it does better only for noise factors above $\eta \approx 0.25$. This is an important property of FLORIS, as measurements in practical deployments are often quite noisy. Running times on the order of 0.1 sec are similar to SLNN's, while SR-LS is much faster (on the order of 2 msec).

\begin{figure}[!h]	
  \centering
  \subfloat[]{%
    \label{RMSE_floris_sim_comparison_2d_add}%
    \includegraphics[width=0.5\linewidth]{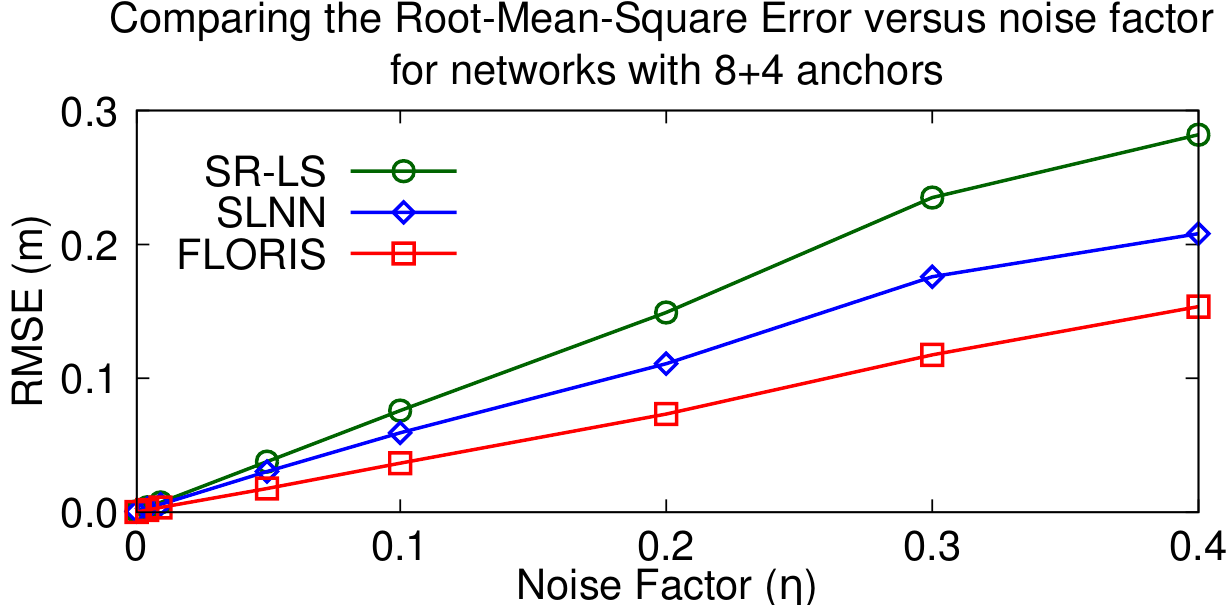}} 
  \subfloat[]{%
    \label{RMSE_floris_sim_comparison_3d_add}%
    \includegraphics[width=0.5\linewidth]{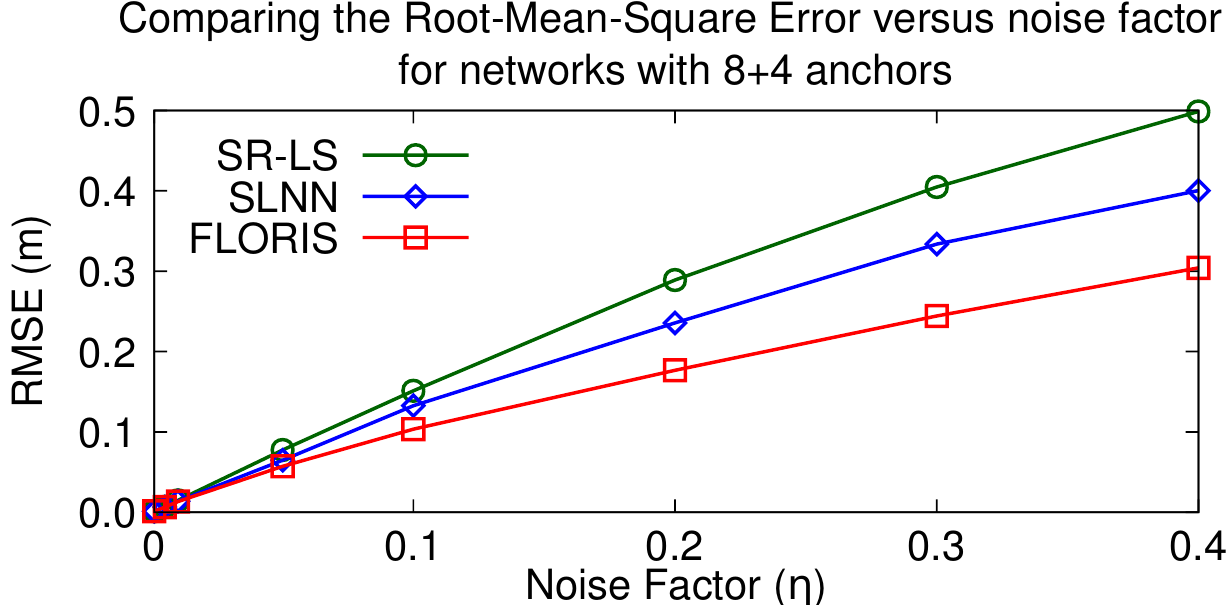}}
  \caption[]{Single-source localization: RMSE vs.\ noise factor for 8 acoustic anchors and 4 additional visual anchors. \subref{RMSE_floris_sim_comparison_2d_add} 2D. \subref{RMSE_floris_sim_comparison_3d_add} 3D.}
  \label{RMSE_floris_sim_comparison_add}
\end{figure}

\noindent\emph{Example 2:} When angular information is derived from a video camera (installed on the source) and fiducial markers (printed on paper and posted on the operating environment), one could argue that the latter contribute negligibly to the total hardware cost of the localization system. Therefore, as long as the computational complexity remains within reasonable limits, it makes sense to investigate how the performance changes when all algorithms use the same fixed number of acoustic anchors, and FLORIS resorts to additional visual anchors. Keeping the same simulation methodology outlined above, Figure \ref{RMSE_floris_sim_comparison_add} shows results for this modified approach with 8 acoustic anchors and 4 additional visual anchors for FLORIS. Given the relatively narrow performance gap between FLORIS and SLNN observed previously in Figure \ref{RMSE_floris_sim_comparison_convert}, it is not surprising to see that the additional observations now give FLORIS a clear advantage, leading to a consistent performance gap to SLNN that is similar to the one between SR-LS and SLNN. The running time of FLORIS is again approximately 0.1 sec.

\noindent\emph{Example 3:} RMSE is the preferred performance metric in our work, but in a semidefinite relaxation method such as FLORIS the quality of the relaxation is also of interest. For the same conditions of Example 2, Table \ref{Rank_1_percentage_diff_8+4networks} lists the relative frequency of instances for the random scenario such that the solution matrix approximately satisfies $\text{rank}(\mathbf{W}) = 1$, in which case the solution of the relaxed problem coincides with that of the original nonconvex problem. We consider matrix $\mathbf{W}$ to effectively have rank 1 when the ratio between its first and second singular values is at least 20.
%2D
\begin{table}[!t]
 \caption{Percentage of Rank-1 FLORIS Solutions for 8+4 Network Configurations vs.\ Noise Factor} \label{Rank_1_percentage_diff_8+4networks}
\footnotesize 
	\centering
	\resizebox{\columnwidth}{!}{%
        \begin{tabular}{ | c | c | c | c | c | c | c | c | c |}
  		\hline
  		\diagbox{Rank $\approx 1$ (\%)}{$\eta$} & 0.001 & 0.005 & 0.01 & 0.05 & 0.1 & 0.2 & 0.3 & 0.4  \\ \hline
  		2D & 100 & 100 & 100 & 99.8 & 95.1 & 89.7 & 77.1 & 80.9 \\ 
  		3D & 99.8 & 99.6 & 99.0 & 91.9 & 83.2 & 77.9 & 76.3 & 69.1 \\ \hline 
	\end{tabular}%
	}
\end{table}
These results show that the likelihood of obtaining a rank 1 solution is high even for large noise factors, confirming the tightness of the relaxation. Problem instances that yield approximately rank 1 solutions lead to accurate source localization, but the degradation is graceful and even solutions that do not meet the rank criterion usually lead to satisfactory localization performance.
\subsubsection{Experimental Results}

An experimental set-up was developed to test FLORIS. This consisted of Cricket~\cite{cricketlocalsys:Priyantha:2000} beacon nodes as acoustic anchors and ARUCO~\cite{Aruco2014} augmented reality tags as visual anchors (see Figure \ref{fig:experimental_setup_glass}). The target, comprising a Cricket listener node rigidly coupled to a video camera, could roam inside a confined volume of about 50 $\mathrm{m}^3$. In addition to the calibration issues described in Section \ref{sec:self-calib}, several practical issues had to be overcome to obtain usable multisensor data, including tuning the directionality of Cricket nodes (using parabolic reflectors to widen the spatial coverage of ultrasonic transducers) and calibrating their range measurements individually. We followed the range error mitigation and position refinement procedure presented in~\cite{Conti:NetworkExperimentation:2012}, which reduced the localization error in an indoor experimental set-up. 

\begin{figure}[!h]	
%\raggedleft
	\centering
	\includegraphics[width=0.6\columnwidth]{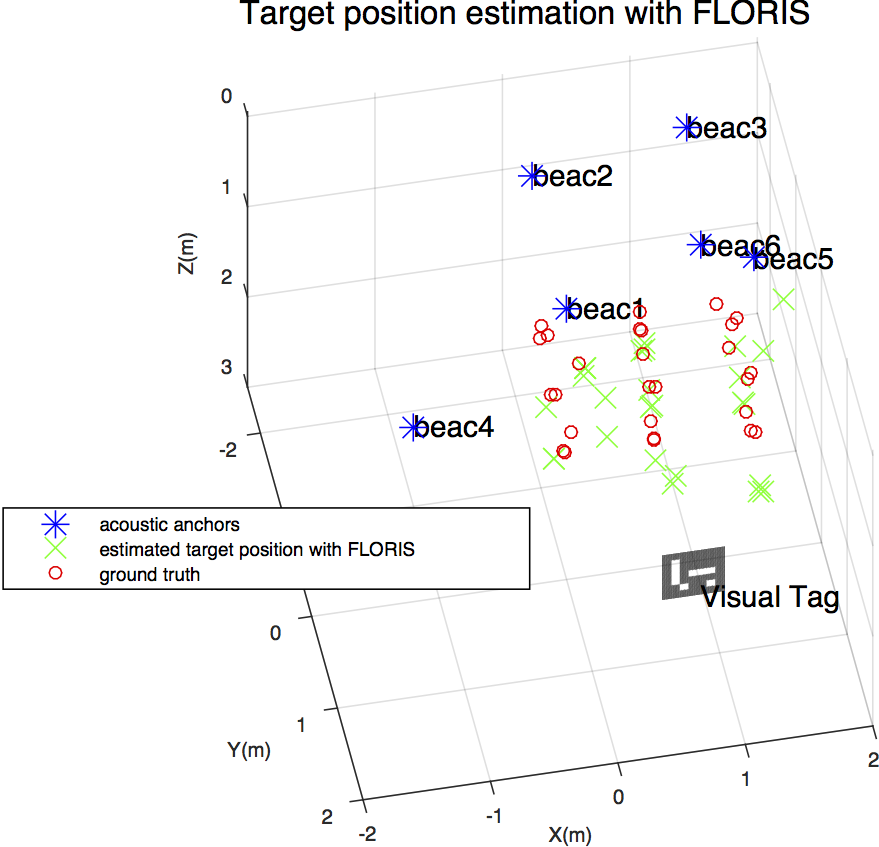}
	\vspace{-2mm}
	\caption{Target position estimates given by FLORIS and ground truth positions (real data). \label{hybrid_target_pos_GT}}
%\vspace{-4mm}
\end{figure}

Several datasets of range and orientation measurements were acquired. Figure \ref{hybrid_target_pos_GT} shows the ground truth and the target positions estimated by FLORIS during a walk through the set-up. It can be observed that the estimated positions are, globally, close to the ground truth\footnote{Due to the limited accuracy of manually measuring target and anchor positions in the room, as well as approximating the Cricket acoustic transmitter/receiver pair as a single point, the ``ground truth'' may actually be affected by errors of 5--10 cm.}.

Table \ref{table:RMSE_localization_cri_calib} (upper row) compares the RMSE obtained for the proposed algorithm with SR-LS and SLNN for a particular dataset comprising 26 target positions. These are consistent with simulation results, showing a moderate advantage of FLORIS over SLNN, and significantly larger errors of SR-LS. Similarly to the simulated scenarios with strong noise, the relaxation used in FLORIS may not be tight mainly due to the poor quality of some of the range estimates. The bottom row of Table \ref{table:RMSE_localization_cri_calib} quantifies the improvement that can be achieved through refinement, i.e., when the computed solution is used as initialization in a numerical minimization procedure that directly operates on the original nonconvex cost \eqref{eq:FLORIS_cost_function}. A similar strategy is used for the other algorithms.
% We performed an extra refinement step of the global solution for all methods using a generic minimization function which attempts to find a local minimizer of the cost function given the previous estimated positions as starting points. The results of this post step are denoted as RMSE-min. 
\begin{table}[!t]
\caption{3D single-source localization performance for an experimental dataset comprising 1 ARUCO (visual) + 6 Cricket (acoustic) anchor nodes.} \label{table:RMSE_localization_cri_calib}
\centering
\begin{tabular}{ c | c | c | c |}   
  \cline{2-4}                   
  &  SR-LS~\cite{sourceLS:beck:2008} & SLNN~\cite{SLCP:pinar:2014} & FLORIS \\ \cline{1-4}
  \multicolumn{1}{ |c| } {RMSE (m)} & 0.43 & 0.33 & 0.31 \\
  \multicolumn{1}{ |c| } {RMSE w/ refinement (m)} & 0.42 & 0.29 & 0.23 \\ \cline{1-4}  
\end{tabular}
%\vspace{-2mm}
\end{table}
% From Table \ref{table:RMSE_localization_cri_calib}, it can be observed that FLORIS provides better accuracy estimates when compared with the other methods, corroborating the good performance already found in simulation.
Refining the solutions of FLORIS and SLNN produces noticeable improvements, but their relative performance gap does not decrease. As for SR-LS, refinement produces virtually no improvement as the algorithm directly addresses its nonconvex cost function without any relaxation.

Due to reliability issues related to the directivity of Cricket nodes, there were regions inside the volume where the number of anchor sightings was insufficient for range-only localization (i.e., fewer than 4 anchors). In this respect FLORIS showed better consistency of coverage, as the detection of visual markers was more robust and could often compensate for a shortage of range measurements at any given target position. We expect this type of advantage to carry over to other operating scenarios and ranging technologies.
%This is a very encouraging result, showing that although several practical limitations were found in the deployment of the experimental set-up (noisy measurements given by the Cricket System and the calibration procedure) they were successfully overcome. 
% These findings support the claim that the hybrid method proposed is a significant contribution to the field, outperforming the previous works in the case of practical noisy data. 

% Overall, we believe these experimental results validate FLORIS as a practically relevant and solid algorithm with appealing accuracy and robustness properties. We stress that FLORIS is potentially more flexible and scalable than other methods operating on a single type of sensed variable.
% More specifically, the hybrid approach is able to successfully localize a target in situations where other algorithms cannot. For example, in 3D, when information from only 3 acoustic anchors is received, it is not possible to estimate a correct position. Yet, if an object (or a tag in our current implementation) is recognized, FLORIS produces a valid estimation.
Furthermore, although FLORIS experiences a well-known deterioration in performance as the source moves outside the convex hull spanned by the anchors, we have observed that the degradation is more progressive and smoother than in range-only algorithms. 
So far, the evidence for this effect remains only empirical, and a more careful characterization should be undertaken.
%Hence, the method proposed paper conveys further flexibility and scalability to localization systems.
\subsection{Cooperative Localization (CLORIS)}

\subsubsection{Numerical Results} 
\begin{figure}[!h]
  \centering
  \subfloat[]{%
    \label{Colab_RMSE_compare_fixed2D}%
    \includegraphics[width=0.5\linewidth]{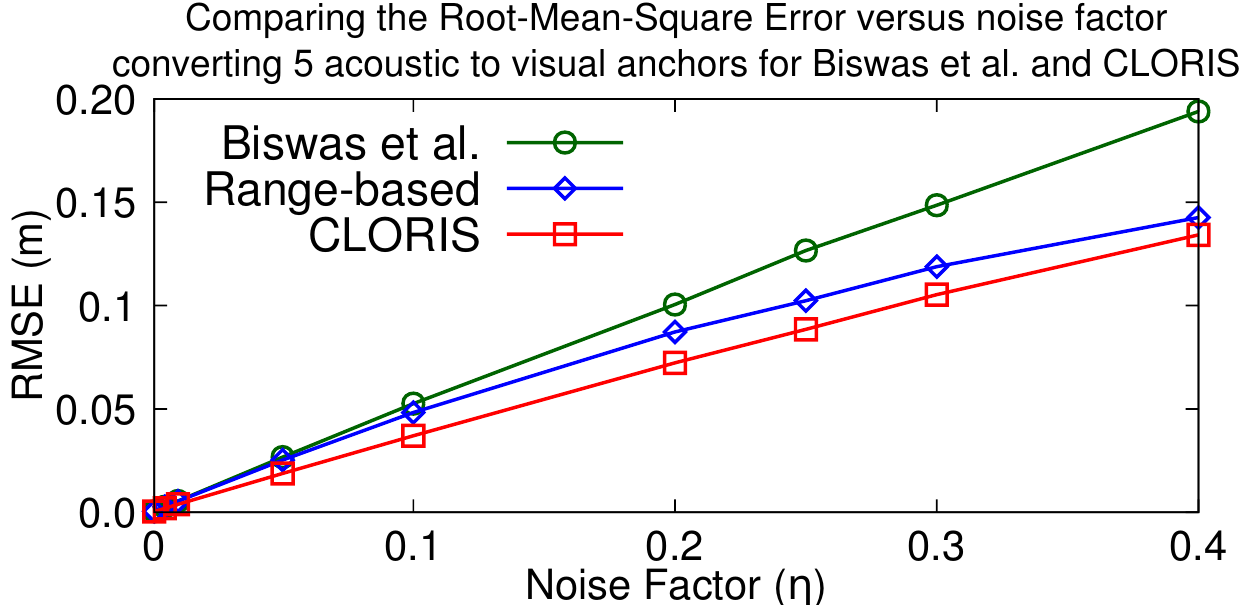}} 
  \subfloat[]{%
    \label{Colab_RMSE_compare_fixed3D}%
    \includegraphics[width=0.5\linewidth]{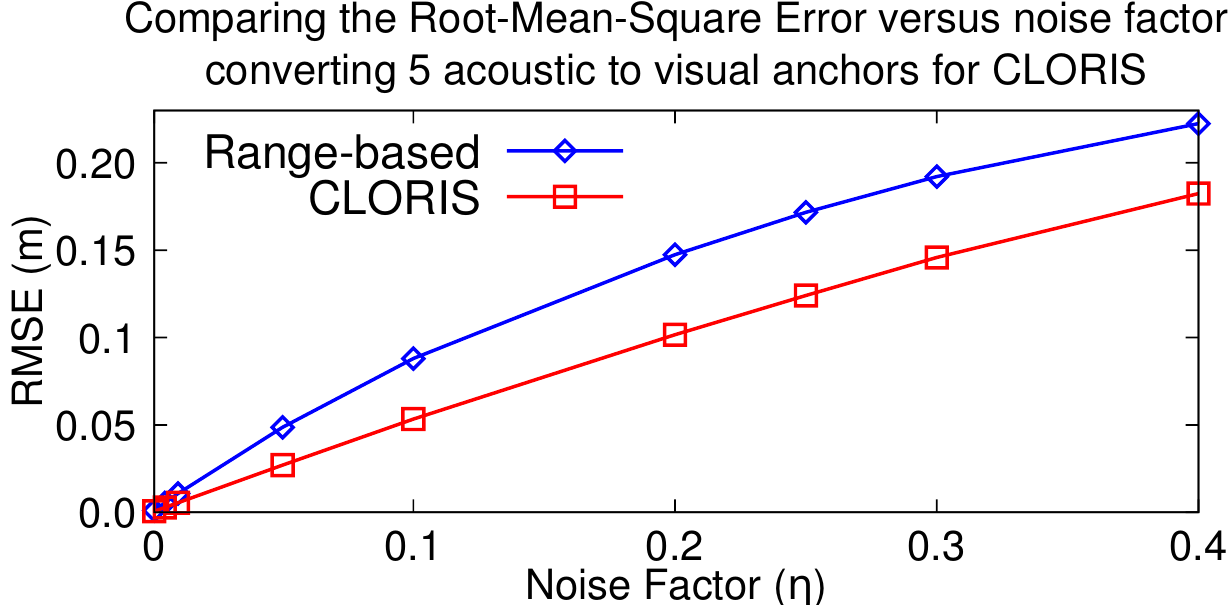}}
	\caption[]{Cooperative localization: RMSE vs.\ noise factor for 4 sensors and 13 acoustic or 8 acoustic + 5 visual anchors. \subref{Colab_RMSE_compare_fixed2D} 2D. \subref{Colab_RMSE_compare_fixed3D} 3D.} 
        \label{Colab_RMSE_compare_fixed}
\end{figure}

~\\ 
\emph{Example 1:} 
Figure \ref{Colab_RMSE_compare_fixed} shows the RMSE versus noise factor for localizable networks comprising 13 acoustic anchors and 4 sensors. Similarly to what was done for the single-source case, in CLORIS and in the method by Biswas et al.~\cite{Biswas05integrationof} 5 anchors were converted into visual ones for a fair comparison with the range-only benchmark algorithm~\cite{soares:coopdistlocal:2014}.  
In 2D, the average error of CLORIS is always lower than the one of~\cite{Biswas05integrationof}, indicating that our relaxation of the original problem is tighter. 
Also, CLORIS outperforms the range-only algorithm for all noise factors both in 2D and 3D, the gap being larger than in FLORIS vs.\ SLNN. This may be partly due to the fact that now averaging occurs over fewer (non-random and localizable) spatial configurations, in which angular information may bring much added value. Both algorithms were stopped when the gradient norm $\| \nabla g_i (w_i) + \nabla h_i(w_i) \|$ reached $10^{-6}$. Running times for CLORIS and for the range-only algorithm were approximately 0.2 sec and 0.5 sec, respectively. The reduction in execution time is a consequence of both the smaller number of iterations performed by CLORIS to meet the stopping criterion and also the lower computational complexity of projecting a point onto a straight line, performed by CLORIS for anchor-sensor/sensor-sensor angle measurements, when compared to the complexity of projecting a point onto a ball, for range measurements.  

\noindent\emph{Example 2:}
In this example we look at the performance gap of CLORIS to non-relaxed (and non-convex) approaches. These have the potential to attain high accuracy, but may converge to undesirable solutions when the underlying cost functions are multimodal. We therefore consider these as refinement (i.e., postprocessing) methods, and initialize them with the sensor positions obtained by CLORIS. Specifically, we benchmark CLORIS against the ML approach of Huang et al.~\cite{Huang:CRLB_AOA:2016}, and a simple algorithm that directly optimizes our non-relaxed hybrid cost function~\eqref{eq:Cost_function_colab_fusion}. The ML formulation of~\cite{Huang:CRLB_AOA:2016} is based on the following proposed likelihood function

%In this example we look at the performance gap of CLORIS to non-relaxed (and non-convex) approaches, and to the Cram\'{e}r-Rao Lower Bound (CRLB) as well. For $N$ sensors we compute the CRLB as $\sqrt{\frac{1}{N}\mathrm{tr}(\mathbf{F}^{-1})}$, where $\mathbf{F}$ is the Fisher Information Matrix (FIM), and use the FIM derived in~\cite{Huang:CRLB_AOA:2016} for 3D collaborative positioning with ranges and bearings. The results of~\cite{Huang:CRLB_AOA:2016} are based on the following proposed log-likelihood function
%\textcolor{red}{In order to further analyze the estimation accuracy we compute the CRLB for the network localization problem.
%We based ourselves in the CRLB derived in~\cite{Huang:CRLB_AOA:2016} for cooperative localization with angle of arrival in conjunction with range measurements.
%The (nonconvex) cost function which the algorithm proposed in~\cite{Huang:CRLB_AOA:2016} iteratively optimizes is given by: 

\begin{equation}
\begin{aligned}
& f(\mathbf{x}) = \sum_{i\stackrel{\mathcal{R}}{\sim} j} \frac{\left(r_{ij} - \|x_i - x_j\|\right)^2}{\sigma_{rij}^2} + \sum_{i} \sum_{k \in \mathcal{R}_i} \frac{\left(r_{ik} - \|x_i - a_k\|\right)^2}{\sigma_{rik}^2} \\
& + \sum_{i\stackrel{\mathcal{T}}{\sim} j} \frac{\left(\phi_{ij} - \text{atan2}\left(x_i(2) - x_j(2),x_i(1) - x_j(1)\right)\right)^2}{\sigma_{\phi ij}^2} \\
& + \sum_{i} \sum_{k \in \mathcal{T}_{i}} \frac{\left(\phi_{ik} - \text{atan2}\left(x_i(2) - a_k(2),x_i(1) - a_k(1)\right)\right)^2}{\sigma_{\phi ik}^2} \\
& + \sum_{i\stackrel{\mathcal{T}}{\sim} j} \frac{\left(\alpha_{ij} - \text{acos}\left(\frac{x_i(3) - x_j(3)}{r_{ij}}\right)\right)^2}{\sigma_{\alpha ij}^2} + \\ & \sum_{i} \sum_{k \in \mathcal{T}_{i}} \frac{\left(\alpha_{ik} - \text{acos}\left(\frac{x_i(3) - a_k(3)}{r_{ik}}\right)\right)^2}{\sigma_{\alpha ik}^2}.
\end{aligned}
  \label{eq:Huang_cost_function}
\end{equation}
In this statistical model the measurements for distances $r_{ij}$, azimuths $\phi_{ij}$ and elevations $\alpha_{ij}$ are independent of each other and across sensor-sensor or sensor-anchor pairs. These are corrupted by zero-mean additive white Gaussian noise with variance $\sigma^2_{\rho ij}$, where $\rho = r$, $\phi$, or $\alpha$ depending on the type of measurement. We note that neither the independence assumption between measured ranges and bearings underlying~\eqref{eq:Huang_cost_function} nor the actual noise statistics match our joint generation model described at the start of Section~\ref{sec:results}. However, since deriving a ML formulation for our model is beyond the scope of this paper, and the mismatch with the assumptions of~\cite{Huang:CRLB_AOA:2016} does not seem severe under weak to moderate noise, we will directly feed that ML estimator with the same synthetic data that CLORIS operates on. To minimize the impact of the mismatch we provide additional side information to the ML estimator in the form of precomputed variances for all individual terms in~\eqref{eq:Huang_cost_function}, empirically obtained over a sizable number of noise realizations for any given geometric configuration.

\begin{figure}[!h]
	\centering
	\includegraphics[width=0.6\linewidth]{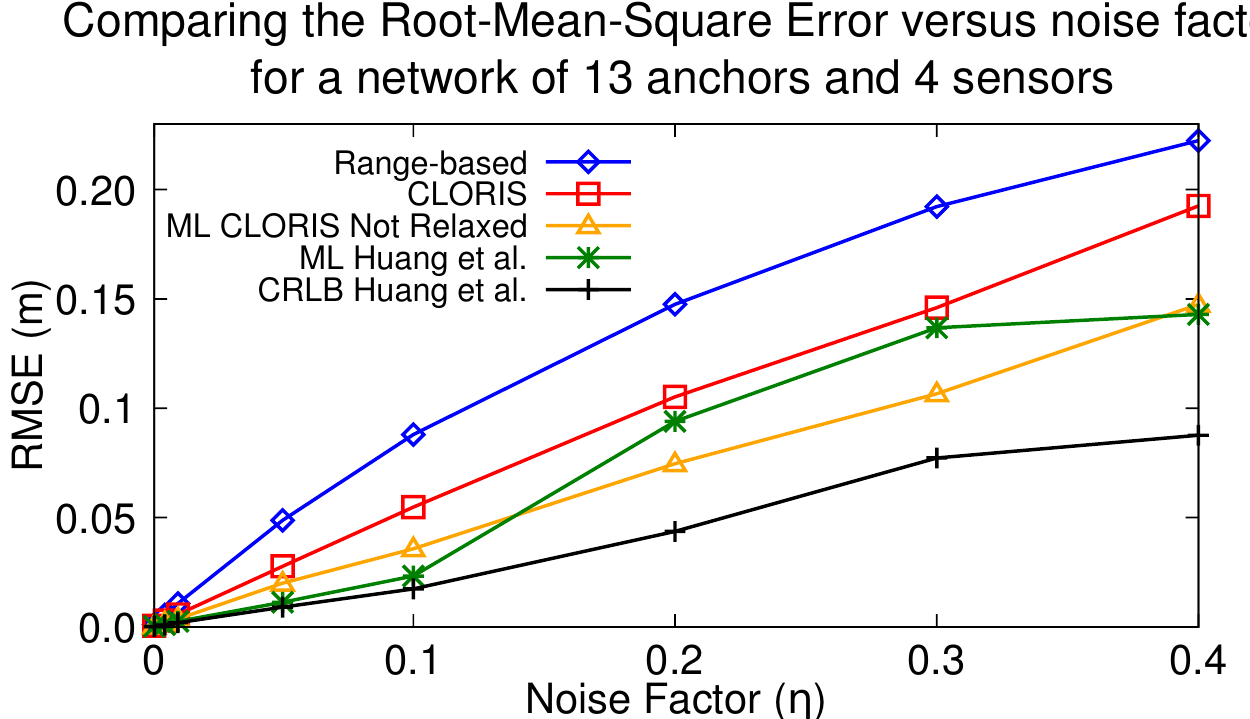}
	\caption{3D Cooperative localization performance: RMSE vs.\ noise factor.} \label{Colab_RMSE_compare_CRLB}
\end{figure}

For the same 3D network configuration of Example 1, Figure \ref{Colab_RMSE_compare_CRLB} plots the RMSE as a function of the noise factor. In addition to CLORIS and the range-only algorithm of~\cite{soares:coopdistlocal:2014}, the figure includes curves pertaining to direct refinement, using one of \emph{MATLAB}'s generic unconstrained nonlinear minimization functions (\texttt{fminsearch}), of our non-relaxed cost function \eqref{eq:Cost_function_colab_fusion} (related to what is developed in~\cite{Piovesan:convex-nonconvex:2016} for range-only) and the likelihood function \eqref{eq:Huang_cost_function} proposed by Huang et al.~\cite{Huang:CRLB_AOA:2016}. To complement the ML results, the figure also shows the Cram\'{e}r-Rao Lower Bound (CRLB) computed as $\sqrt{\frac{1}{N}\mathrm{tr}(\mathbf{F}^{-1})}$, where $\mathbf{F}$ is the Fisher Information Matrix (FIM) given in~\cite{Huang:CRLB_AOA:2016} for 3D collaborative positioning with ranges and bearings. Again, due to model mismatch there are no guarantees that this actually represents the best attainable unbiased performance with our simulated data, but for moderate noise it does provide a useful order-of-magnitude notion of how inefficient the ML estimator might be in the non-asymptotic regime.

%For the same 3D network configuration of Example 1, Figure \ref{Colab_RMSE_compare_CRLB} plots the RMSE as a function of the noise factor for CLORIS, the range-only algorithm of~\cite{soares:coopdistlocal:2014}, and the CRLB obtained as described above. The figure also includes curves pertaining to direct refinement, using one of \emph{MATLAB}'s generic unconstrained nonlinear minimization functions (\texttt{fminsearch}), of our non-relaxed cost function \eqref{eq:Cost_function_colab_fusion} and the one proposed by Huang et al. \eqref{eq:Huang_cost_function}.

The results show similar performance for ML-based estimation using \eqref{eq:Huang_cost_function} and refinement using our non-relaxed cost function \eqref{eq:Cost_function_colab_fusion}. This suggests that, at least for some relevant scenarios, our ad-hoc quadratic hybrid cost function is a reasonable surrogate for a ``true'' likelihood, while having better analytical tractability and more parsimonious parameterization. The narrow performance gap between CLORIS and both refinement schemes confirms the effectiveness and accuracy of our convex relaxation strategy. Finally, although the curves for ML and the CRLB should be interpreted with caution for the reasons outlined above, we note that their gap is relatively narrow, suggesting that drastic performance improvement over the former with alternative estimation approaches is not expectable.

In light of the simulation results, and to the extent possible in the absence of a CRLB for our data model, we argue that it seems unlikely that the actual performance lower bound for unbiased estimation will be drastically lower than the CRLB shown in Figure \ref{Colab_RMSE_compare_CRLB}. Our conclusion is that the gap between CLORIS and a true ML estimator, even if it were efficient, would remain narrow.

%The results show a gap from CLORIS to the CRLB, but not an overly wide one. Importantly, the ML estimate based on the likelihood function of Huang et al.~\cite{Huang:CRLB_AOA:2016} does not touch the CRLB for this problem size, either, and in fact only improves moderately on the initialization provided by CLORIS. Also, direct optimization of our non-relaxed cost \eqref{eq:Cost_function_colab_fusion} yields similar performance to that obtained from the likelihood function \eqref{eq:Huang_cost_function}. This indicates that, for reasonable scenarios, our ad-hoc hybrid cost function is a reasonable surrogate for the likelihood, while having better analytical tractability and more parsimonious parameterization.

%\textcolor{red}{Figure \ref{Colab_RMSE_compare_CRLB} plots, together with the RMSE of CLORIS and the range-based method, both the CRLB and the RMSE of the ML function \eqref{eq:Huang_cost_function} from~\cite{Huang:CRLB_AOA:2016}, as well as the RMSE of the original (non relaxed) and the convex (relaxed) versions of the CLORIS cost function, \eqref{eq:Cost_function_colab_fusion} and \eqref{eq:Cost_function_colab_fusion_final}, respectively, for the same 3D network configuration of \textit{Example 2}.
%Figure \ref{Colab_RMSE_compare_CRLB} shows that not only the original CLORIS ML function is relatively close to the CRLB, but also the gap induced by the performed relaxation is not large. Furthermore, CLORIS method achieved a performance extremely close to the relaxed ML function.}

\begin{figure}[!h]
	\centering
	\includegraphics[width=0.6\linewidth]{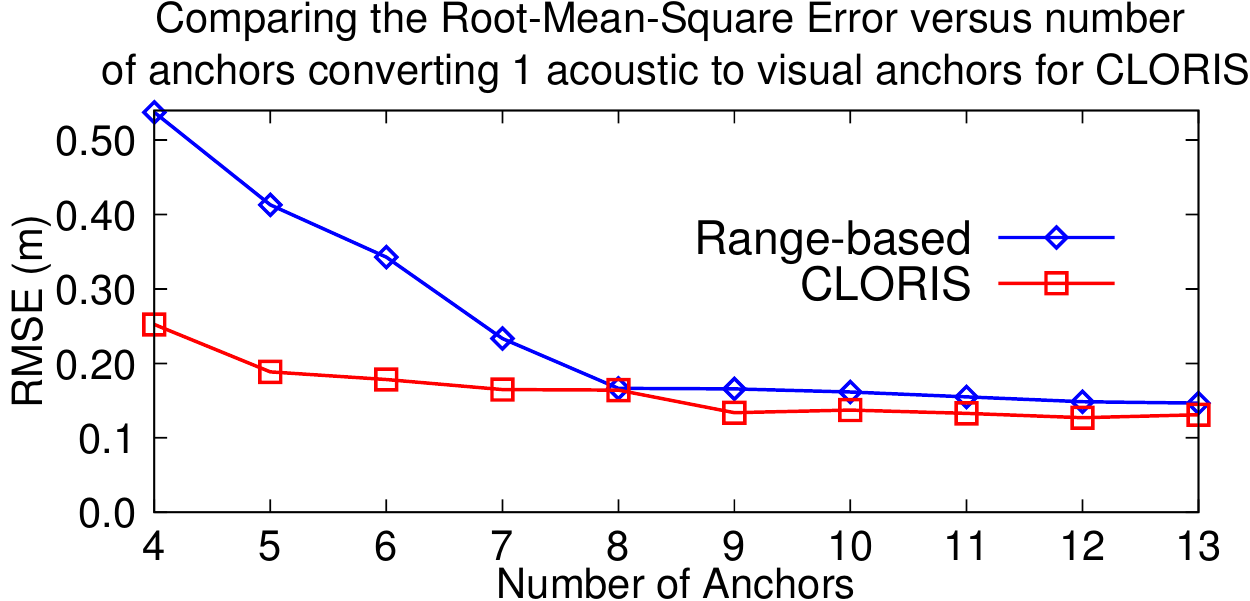}
	\caption{Cooperative 3D localization: RMSE as a function of the number of anchors for 4 sensors and noise factor $\eta = 0.2$. For CLORIS, one acoustic anchor is converted to visual.} \label{Colab_RMSE_compare_fixed3D_numanchors}
\end{figure}

\noindent\emph{Example 3:}
Figure \ref{Colab_RMSE_compare_fixed3D_numanchors} shows the performance of CLORIS and the range-only benchmark in 3D as a function of the number of anchors, including 1 visual anchor for CLORIS, 4 sensors and noise factor $\eta = 0.2$.
As expected, increasing the number of anchors gradually improves the accuracy of both algorithms, but CLORIS retains an advantage for all network sizes. For smaller networks (up to 7 anchors), in particular, the results show that including angular information even for a single anchor is very beneficial, as it tends to provide a good fix to the positions of the sensors that observe it with few additional range measurements. This provides a compelling argument in favor of hybrid localization methods that can make the best of what measurement resources are available in a given configuration.
 % With more visual anchors (not shown) the solution obtained by CLORIS is defined mostly by the angular information, and improves only slightly as the number of acoustic anchors increases.
 
\noindent\emph{Example 4:} Since CLORIS may be used for the specific case of networks with only one sensor (single-source), the question arises as to whether it outperforms FLORIS, their computational complexities being comparable, despite very different algorithmic structures. Solving the single-source localization problem for randomly generated configurations of 8 acoustic and 4 visual anchors produced the results in Table \ref{SS_Perform_FLORIS_CLORIS}. In this scenario, FLORIS consistently shows better performance in 2D, and outperforms CLORIS in 3D for larger noise factors. We therefore argue that both algorithms are relevant for their targeted localization scenarios.
\begin{table}[!t]
 \caption{Single-source localization performance comparison between FLORIS and CLORIS.} \label{SS_Perform_FLORIS_CLORIS}
\footnotesize 
	\centering
	\resizebox{\columnwidth}{!}{%
          \begin{tabular}{ | c | c | c | c | c | c | c | c |}
            \hline                       
            \multicolumn{2}{|c|}{\diagbox{RMSE}{$\eta$}} & 0.001 & 0.005 & 0.01 & 0.05 & 0.1 & 0.2 \\ \hline
            \multirow{2}{*}{2D}
             & FLORIS & 0.0004 & 0.0019 & 0.0038 & 0.0193 & 0.0374 & 0.0755 \\
             % & SLNN & 0.0004 & 0.0022 & 0.0044 & 0.0224 & 0.0413 & 0.0802 \\ \cline{2-8}
             & CLORIS & 0.0008 & 0.0044 & 0.0088 & 0.0443 & 0.0861 & 0.1451 \\ \hline\hline 
            \multirow{2}{*}{3D}
  	    & FLORIS & 0.0014 & 0.0066 & 0.0126 & 0.0568 & 0.1032 & 0.1814 \\
            % & SLNN & 0.0009 & 0.0045 & 0.0088 & 0.0444 & 0.0850 & 0.1659 \\ \cline{2-8}
            & CLORIS & 0.0011 & 0.0053 & 0.0108 & 0.0527 & 0.1036 & 0.1984 \\ \hline 
	\end{tabular}%
	}
\end{table}

%Firstly, we compare the performance of FLORIS and CLORIS for the single source scenario. We show that the first attains lower error, thus justifying/supporting, for this case, the choice of FLORIS framework (relaxation of the cost function to an SDP) over the disk-relaxation approach.

\subsubsection{Experimental Results}
For testing CLORIS experimentally, we deployed two targets (each comprising, similarly to the single-source scenario, a listener Cricket node mounted on a camera, as well as fiducial markers placed on the outer surface of the enclosure) in the area covered by the acoustic and visual anchors. We collected a dataset of selected sensor-sensor and sensor-anchor range and angle measurements. Due to limitations of Cricket nodes, whose uncostumized firmware does not allow switching between transmit (beacon) and receive (listener) modes, internode range measurements are not available, which excludes true cooperation in range-only localization. The results reported here actually focus on particular configurations where the network is localizable using all measurements, but not if only range measurements are considered.

Table \ref{table:RMSE_networklocalization_mixed} shows the RMSE for range-only~\cite{soares:coopdistlocal:2014} and hybrid (CLORIS) localization. Because we are interested in assessing the impact of angular measurements as the errors propagate through rotations and translations to the global reference frame (see Section \ref{sec:self-calib}), the table also shows the RMSE computed using half-real data: actual Cricket range measurements but synthetic angular measurements derived from the ground truth.
\begin{figure}[!h]	
	\centering
	\includegraphics[width=0.65\columnwidth]{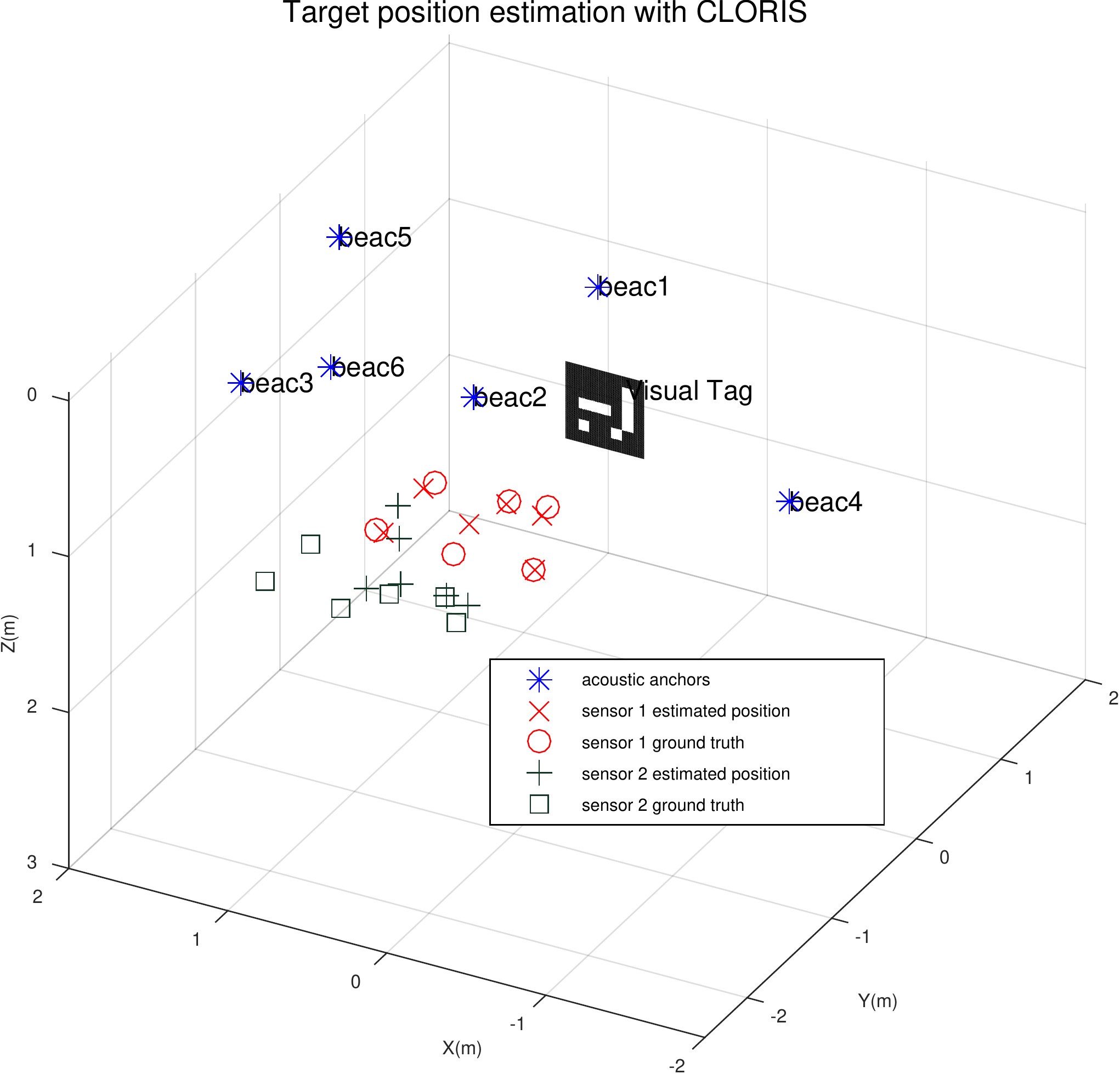}
	\caption{CLORIS estimated sensor positions vs.\ ground truth positions during a walk through the experimental set-up. \label{Colab_estpos_final}}
\end{figure}  
\begin{table}[!t]
\caption{3D network localization performance for an experimental dataset comprising 2 sensors and 1 ARUCO (visual) + 6 Cricket (acoustic) anchor nodes.}\label{table:RMSE_networklocalization_mixed}
\centering
\begin{tabular}{c|c|ccc}
\cline{2-2}
                                                                 & RMSE (m) &  &  &  \\ \cline{1-2}
\multicolumn{1}{|c|}{Range-based Disk Relaxation}                   & 1.43  &  &  &  \\
\multicolumn{1}{|c|}{CLORIS} & 0.55       &  &  &  \\
\multicolumn{1}{|c|}{CLORIS w/ synthetic angles}    & 0.50  &  &  &  \\ \cline{1-2}
\end{tabular}
\end{table}
Compared with the results for single-source localization in Table \ref{table:RMSE_localization_cri_calib}, cooperation in CLORIS is seen to provide enhanced precision. Note that the table shows the \emph{total} RMSE \eqref{eq:rmse}, and not the RMSE per sensor. Large errors occur in range-only localization as the coordinates of one of the sensors cannot be univocally reconstructed from the available ranges in some of the snapshots. This eloquently underscores the potential benefits of hybrid localization schemes in applications. Figure \ref{Colab_estpos_final} shows the sensor positions computed by CLORIS in the dataset, in good agreement with the ground truth.

The degradation of CLORIS with respect to the variant using synthetic angles is surprisingly modest, bearing in mind the fact that the sensor-sensor direction is computed in the global reference frame through the concatenation of two sets of (noisy) rotations/translations, namely, sensor 2 to sensor 1, followed by sensor 1 to the visual anchor. Still, at the individual sensor level we note that there is a significant discrepancy between the localization performance for sensor 1, which directly observes a visual anchor, and sensor 2, for which the only available angular information is derived through the above two-step procedure and therefore has lower quality. These differences are clearly visible in Figure \ref{Colab_estpos_final}, where the estimated positions of sensor 1 are closer to the ground truth than those of sensor 2. In particular, the accuracy for sensor 1 is higher than in single-source localization using FLORIS with the same set of anchors.

\section{Conclusions and Future Work}
\label{sec:conclusions}

We addressed the problem of indoor localization using ranges and angles through a least-squares approach that seamlessly integrates these heterogeneous measurements into a single optimization problem. While our (nonlinear) quadratic cost is not a likelihood function, thus suboptimal from the perspective of estimation theory, it is mathematically tractable and was shown to lead to state-of-the-art accuracy in simulation and in real experiments. We have addressed both the single-source and network (cooperative) localization paradigms, and developed efficient algorithms for each of them (both extendable to arbitrary space dimensions). These are genuine fusion algorithms that estimate node positions from the \emph{full} set of measurements, thus bypassing a critical assumption in other published hybrid algorithms that require the nodes to be localizable using only measurements of either type. However, we currently do impose range-only localizability for the optional (one-off) automated calibration procedure that is used to set up the conversion of angular measurements from camera-centric coordinates to a global reference frame. Eliminating this constraint might be a topic for future developments.

FLORIS, the algorithm for single-source localization based on semidefinite relaxation, was shown to outperform benchmark methods in simulated 2D scenarios, and also in 3D under medium/strong measurement noise. Experimental results in our testbed validated FLORIS as a practically relevant and solid algorithm with appealing accuracy and robustness properties.

We stress that FLORIS is potentially more flexible and scalable than other methods operating on a single type of sensed variable. This was quite evident in our testbed, where limitations in the technology used for acoustic ranging made it somewhat cumbersome to obtain enough measurements at arbitrary source positions within the area of coverage for 3D range-only localization, whereas the presence of a single visual anchor simplified things considerably for hybrid localization. Surely, this does not exclude the existence of topologies in which including visual anchors may deteriorate accuracy. Nonetheless, throughout our experiments we have never encountered such cases. On the contrary, we found that using bearing measurements consistently reduces the typical accuracy degradation that occurs when the source moves outside the convex hull spanned by the anchors. 

CLORIS, the algorithm for network localization based on a so-called disk relaxation and accelerated Nesterov's gradient descent, similarly improved upon the performance of its range-only counterpart (which itself was previously shown to beat several benchmark algorithms) in both 2D and 3D. Experimental results in our testbed confirmed that cooperation adds value in practice, reducing the error in computed node positions relative to the single-source case. 

CLORIS has a highly parallelizable computational structure that inherently makes it scalable as the size of the network increases. In this work we have not explored the related issue (definitely worth pursuing) of creating a \emph{distributed} version of CLORIS, that would also involve local processing of measurements to express the spatial information in the global reference frame. The issue may be trivial for other sensing modalities and/or using technological aides such as IMUs.

We tested CLORIS for single-source localization with good results, although the precision was lower than that of FLORIS in most simulated conditions. Running times were moderate but larger for CLORIS (about 0.2 sec in our experimental data vs.\ 0.1 sec for FLORIS), although the situation might be reversed for larger problems, as the complexity of semidefinite relaxation problems increases faster with the number of nodes.   

Currently, angular and range terms are weighted identically on the hybrid least-squares cost function. Introducing different weights, depending on the relative precisions of range and angular sensors, seems quite doable and might afford practical gains. Along the same vein, developing and assimilating into the cost criterion a more realistic model for angular measurements obtained from fiducial markers might prove beneficial, as the ``sweet spot'' for these is known to occur at slant angles in the vicinity of 45 degrees. Finally, a longer-term goal would be to dispense with fiducial markers altogether, and derive angular information through opportunistic detection/recognition of objects or structures present in the environment by computer vision techniques. This would subsume, but considerably expand, the proposed self-calibration formulation for assimilation of range and video data.

\section*{Acknowledgements}
This research was partially supported by Funda\c{c}\~{a}o para a Ci\^{e}ncia e a Tecnologia (project UID/EEA/50009/2013) and EU FP7 project WiMUST (grant agreement no.\ 645141).

%\section*{References}

\bibliography{journal_paper}

\end{document}